\documentclass[12pt]{amsart}
\usepackage{amssymb}
\usepackage{verbatim}
\usepackage{euscript}

\let\cal\mathcal

\makeatletter \@mparswitchfalse \makeatother

\newtheorem{theorem}{Theorem}
\newtheorem{lemma}[theorem]{Lemma}
\newtheorem{corollary}[theorem]{Corollary}
\newtheorem{proposition}[theorem]{Proposition}
\theoremstyle{remark}
\newtheorem{remark}[theorem]{Remark}
\newtheorem{remarks}[theorem]{Remarks}
\theoremstyle{definition}
\newtheorem{definition}[theorem]{Definition}
\newtheorem{problem}[theorem]{Problem}

\numberwithin{equation}{section}
\numberwithin{theorem}{section}

\begin{document}

\title[martingale transforms]{Non-commutative martingale transforms}
\author[N. Randrianantoanina]{Narcisse Randrianantoanina}
\address{Department of Mathematics and Statistics, Miami University, Oxford,
Ohio 45056}
\thanks{Supported in part by NSF grant DMS-0096696 and by a Miami University
Summer Research Appointment.} \email{randrin@muohio.edu}

\def\D{\mathbb D}
\def\S{\mathbb S}
\def\N{\mathbb N}
\def\C{\mathbb C}
\def\J{\mathbb J}
\def\Z{\mathbb Z}
\def\B{\mathbb B}
\def\R{\mathbb R}
\def\M{\cal{M}}
\def\H{\cal{H}}
\let\<\langle
\def\ch{\raise 0.5ex \hbox{$\chi$}}
\def\T{\tau}
\def\E{\cal{E}}

\let\\\cr
\let\phi\varphi
\let\union\bigcup
\let\inter\bigcap
\let\epsilon\varepsilon
\def\supp{\operatorname{supp}}
\def\Im{\operatorname{Im}}
\def\dim{\operatorname{dim}}
\def\Span{\operatorname{span}}
\def\cord{\operatorname{cord}}
\def\Re{\operatorname{Re}}
\def\sqn{\operatorname{sqn}}
\def\log{\operatorname{log}}

\begin{abstract} We prove that non-commutative martingale transforms are
of weak type $(1,1)$.  More precisely, there is an absolute
constant $C$ such that  if $\M$ is a semi-finite von Neumann
algebra and $(\M_n)_{n=1}^\infty$ is an increasing filtration of
von Neumann subalgebras of $\M$ then for any non-commutative
martingale $x=(x_n)_{n=1}^\infty$ in $L^1(\M)$, adapted to
$(\M_n)_{n=1}^\infty$, and any sequence of signs
$(\epsilon_n)_{n=1}^\infty$,
$$\left\Vert \epsilon_1 x_1 + \sum_{n=2}^N \epsilon_n(x_n -x_{n-1})
\right\Vert_{1,\infty} \leq C \left\Vert x_N \right\Vert_1 $$ for
every $N\geq 2$. This generalizes a result of Burkholder from
classical martingale theory to non-commutative setting and answers
positively a question of Pisier and Xu. As applications, we get
the optimal order of the UMD-constants of the Schatten class $S^p$
when $p \to \infty$. Similarly, we prove that the UMD-constant of
the finite dimensional Schatten class $S_n^{1}$ is of order
$\log(n+1)$. We also discuss the Pisier-Xu non-commutative
Burkholder-Gundy inequalities.
\end{abstract}

\maketitle

\section{Introduction}
Non-commutative (or quantum) probability has developed into an
independent field of mathematical research and has received
considerable progress in recent years.  We refer to the books
\cite{AvW} and \cite{Mey} for connections between mathematical
physics, non-commutative probability and classical probability,
the books of Voiculescu, Dykema and  Nica \cite{VDN} and Hiai and
Petz \cite{HiPe} for interplay between operator algebras and free
probability theory, the work of Biane and Speicher \cite{BS}  on
stochastic analysis and free Brownian motion.

In this paper, our main interest is on non-commutative
martingales. Non-commutative martingales have been studied by
several authors.  For instance, pointwise convergence of
non-commutative  martingales was considered in \cite{Cuc} and
\cite{DN}. In \cite{PX}, Pisier and Xu proved a non-commutative
analogue of the Burkholder-Gundy square function inequalities.
Shortly after, Pisier \cite{PS5}, using combinatorial method,
extended their result to a more general class of sequences called
$p$-orthogonal sums when $p$ is an even integer.  Very recently,
Junge and Xu \cite{JX} considered the non-tracial case of the main
result of \cite{PX} along with several related inequalities such
as non-commutative analogue of the classical Burkholder
inequalities on the conditioned square functions among others.
Junge proved in \cite{Ju} non-commutative versions of Doob's
maximal inequalities. We remark that most inequalities considered
in the aforementioned papers were for $p>1$. We continue this line
of research by studying  martingale transforms of non-commutative
bounded $L^1$-martingales.  In the classical probability, the
theory of martingale transforms is well-established and has been
proven to be a very powerful tool not only in probabilistic
situations but also in several parts of analysis.  We refer to the
survey \cite{Bu1} for discussions on this classical topic.  For
instance, Burkholder  \cite{Bu4} proved  that classical martingale
transforms are of weak type (1,1).  Our main result (see
Theorem~\ref{main} below) is a non-commutative analogue of this
classical fact: non-commutative martingale transforms are bounded
as maps from non-commutative $L^1$-spaces into the corresponding
non-commutative weak-$L^1$-spaces. We should point out that this
question was explicitly raised  by Pisier and Xu in the recent
survey \cite{PX3} (Problem 7.5) as it is closely related to the
main result of \cite{PX}. Indeed, combined with general theory of
interpolations of operators of weak types, our main result implies
that for $p>1$, martingale difference sequences in non-commutative
$L^{p}$-spaces are unconditional which in turn imply the
non-commutative Burkholder-Gundy inequalities. This alternative
approach yields constants which are $O(p)$ when $p \to \infty$.
This is explained in Sect.~5. Another application of the main
result is on UMD-constants of non-commutative $L^p$-spaces. It is
now a well known fact that non-commutative $L^p$-spaces on
semi-finite von Neumann algebras are UMD-spaces. The UMD-constants
of these spaces recorded in the literature thus far seems to be of
order $O(p^2)$ when $p \to \infty$. Using the estimates on the
constant of unconditionality of non-commutative martingale
difference sequences, we can deduce that the UMD-constants for
non-commutative $L^p$-spaces are of order $O(p)$ when $p \to
\infty$. We refer to Sect.~4 below for more discussion on this
along with some related results.

The study of martingales in non-commutative cases  often requires
additional insights. In fact, most of usual techniques used in the
classical case are relaying on stopping times or some other basic
truncations which, in many situations, are not available for the
non-commutative setting. Our proof is completely self-contained.
It is  based on a maximal inequality type result from a paper of
Cuculescu \cite{Cuc} (see Proposition~\ref{maximal} below) which
allows ones to reduce the case of  bounded $L^1$-martingales to
bounded $L^2$-supermartingales. Although, such reduction to
supermartingales is standard in classical martingale theory (see
for instance \cite[Chap.~5]{EG}),  the non-commutative setting
presents considerable additional technical difficulty and
therefore requires special care.

The paper is organized as follows:  in Sect.~2 below, we set some
basic preliminary background concerning  non-commutative spaces
and martingale theory that will be needed throughout. Sect.~3 is
devoted mainly to the statement and proof of the main result. In
Sect.~4, we discuss the UMD-constants of non-commutative spaces.
As mentioned above, we revisit the non-commutative Burkholder
inequalities with special attention given to the order of growths
of the constants involved in Sect.~5 and in the last section, we
discuss the class $L\log L$ and formulate some related open
questions.

Our notation and terminology are standard as may be found in the
books \cite{KR1} and \cite{TAK}.

\section{Preliminaries}
Let $\M$ be a semi-finite von Neumann algebra with a normal
faithful semi-finite trace  $\T$. For $1\leq p\leq \infty$, let
$L^p(\M,\T)$ be the associated non-commutative $L^p$-space. Note
that if $p=\infty$, $L^\infty(\M,\tau)$ is just $\M$ with the
usual operator norm; also recall that for $1\leq p<\infty$, the
norm on $L^p(\M,\T)$ is defined by
$$\Vert x \Vert_p =( \T(|x|^p))^{1/p}, \qquad x \in L^p(\M,\T),$$
where $|x|=(x^*x)^{1/2}$ is the usual modulus of $x$.

In order to describe all the spaces involved in this paper, we
recall the general construction of non-commutative spaces as sets
of densely defined operators on a Hilbert space. Throughout, $H$
will denote a Hilbert space and $\M \subseteq B(H)$. The identity
element of $\M$ is denoted by ${\bf 1}$.
 A closed densely defined operator $a$ on $H$ is said to be {\em
affiliated with} $\M$ if $u^* au = a$ for all unitary $u$ in the
commutant $\M'$ of $\M$. If $a$ is a densely defined self-adjoint
operator on $H$, and if $a = \int^\infty_{- \infty} s d e^a_s$ is
its spectral decomposition, then for any Borel subset $B \subseteq
\R$, we denote by $\chi_B(a)$ the corresponding spectral
projection $\int^\infty_{- \infty} \chi_B(s) d e^a_s$. A closed
densely defined operator $a$ on $H$ affiliated with $\M$ is said
to be {\em $\T$-measurable} if there exists a number $s \geq 0$
such that $\T(\chi_{(s, \infty)} (|a|)) < \infty$.

The set of all $\T$-measurable operators will be denoted by
$\overline{\M}$. The set $\overline{\M}$ is a $*$\!-algebra with respect to the
strong sum, the strong product, and the adjoint operation \cite{N}. %
For $x \in \overline{\M}$, the generalized singular value function $\mu (x)$
of $x$ is defined by
$$
\mu_t(x) = \inf \{ s \geq 0: \T(\chi_{(s, \infty)} (|x|)) \leq t \},
\quad \text{ for } t \geq 0.
$$
The function $t \to \mu_t(x)$ from $(0, \T({\bf 1}))$ to $[0, \infty)$
is right continuous, non-increasing and is the inverse of the
distribution function $\lambda (x)$, where $\lambda_s(x) = \T(\chi_{(s,
\infty)}(|x|))$, for $s \geq 0$. For a complete study of $\mu(.)$ and
$\lambda(.)$,  we refer to \cite{FK}.
 For the definition below, we refer
the reader to \cite{BENSHA} and \cite{LT} for the theory of
rearrangement invariant function spaces.

\begin{definition}
Let $E$ be a rearrangement invariant (quasi-) Banach function
space on $(0, \T({\bf 1}))$. We define the symmetric space $E(\M,
\T)$ of measurable operators by setting:
\begin{align*}
E(\M, \T) &= \{ x \in
\overline{\M}\ : \ \mu(x) \in E \} \quad \text{and} \\
\|x\|_{E(\M,\T)} &= \| \mu(x)\|_E,\  \text{ for } x \in E(\M,\T).
\end{align*}
\end{definition}

It is well known that $E(\M, \T)$ is a Banach space (resp.
quasi-Banach space) if $E$ is a Banach space (resp. quasi-Banach
space). The space $E(\M, \T)$ is often referred to as the
non-commutative analogue of the function space $E$ and
  if $E = L^p(0, \T({\bf 1}))$, for $0 < p \leq \infty$, then
$E(\M, \T)$ coincides with the usual non-commutative $L^p$\!-space
associated with $(\M, \T)$. We refer to \cite{CS}, \cite{DDP1},
\cite{DDP3} and \cite{X}
 for more detailed discussions about these
spaces. Of special interest in this paper is the non-commutative
weak $L^1$-space, denoted by $L^{1,\infty}(\M,\T)$ which is
defined as the linear subspace of all $x \in \overline{\M}$ for
which  the quasi-norm
$$ \Vert x \Vert_{1,\infty} := \sup_{t>0} t\mu_t(x) =
\sup_{\lambda>0}\lambda \T(\ch_{(\lambda,\infty)}(|x|)) $$
is finite. Equipped with the quasi-norm $\Vert \cdot \Vert_{1,\infty}$,
$L^{1,\infty}(\M,\T)$ is a quasi-Banach space  and
$\Vert x \Vert_{1,\infty} \leq \Vert x\Vert_1 $ for all $x \in L^1(\M,\T)$.

We now recall the general setup for martingales. The reader is
referred to \cite{Doob} and \cite{GA} for the classical martingale
theory. Let $(\M_n)_{n=1}^\infty$ be an increasing  sequence of
von Neumann subalgebras of $\M$ such that the union of $\M_n$'s is
weak$^*$-dense in $\M$. For each $n\geq 1$, assume that there  is
a conditional expectation $\E_n$
 from $\M$ onto $\M_n$ satisfying:
 \begin{itemize}
 \item[(i)] $\E_n(axb)=a\E_n(x)b$ for all $a, b \in \M_n$ and $x
 \in \M$;
 \item[(ii)] $\T\circ\E_n = \T$.
 \end{itemize}
  It is clear that for every $m$ and  $n$ in $\N$,
$\E_m\E_n=\E_n\E_m =\E_{\min(n,m)}$. Since $\E_n$ is trace
preserving,  it extends to a contractive projection from
$L^p(\M,\T)$ onto $L^p(\M_n,\T_n)$ for all $1\leq p \leq \infty$
where $\T_n$ is the restriction of $\T$ on $\M_n$.
 More generally, a
simple interpolation argument would prove that if $E$ is a rearrangement invariant
Banach function space on $(0, \T({\bf 1}))$ then $\E_n$ is a contraction from
$E(\M,\T)$ onto $E(\M_n,\T_n)$.

Remark that if $\M$ is finite, such conditional expectations
always exist. Indeed, if
  $\cal{N}$ is a von Neumann subalgebra of $\M$. The
embedding $\iota: L^1(\cal{N},\T) \to L^1(\M,\T)$ is an isometry
and the dual map $\E=\iota^*: \M \to \cal{N}$ yields a conditional
expectation (see for instance, \cite[Theorem~3.4]{TAK}).
\begin{definition}
A non-commutative martingale with respect to the filtration $(\M_n)_{n=1}^\infty$ is a
sequence $x=(x_n)_{n=1}^\infty$ in $L^1(\M,\T)$ such that:
$$\E_n(x_{n+1})=x_n \qquad \text{for all}\ n\geq 1.$$
\end{definition}
If additionally $x \in L^p(\M,\T)$ then $x$ is called a $L^p$-martingale.
In this case, we set
$$\Vert x \Vert_p=\sup_{n\geq 1} \Vert x_n \Vert_p.$$
If $\Vert x \Vert_p<\infty$, then $x$ is called a bounded $L^p$-martingale.
The difference sequence of  a martingale $x$ is defined as $dx=(dx_n)_{n=1}^\infty$ with
$dx_1=x_1$ and $dx_n = x_n - x_{n-1}$ for $n\geq 2$.

\medskip

Recall that a subset $K$ of $L^1(\M,\T)$ is said to be uniformly integrable if it is bounded
and for every sequence of projections $(p_n)_{n=1}^\infty $ with
$p_n \downarrow_n 0$, we have
$\lim_{n \to \infty } \sup\{\Vert p_n h p_n \Vert_1; h \in K \}=0$. It is clear that
a  martingale $x=(x_n)_{n=1}^\infty$ in $L^1(\M,\T)$ is uniformly integrable
if and only if there exists $x_\infty \in L^1(\M,\T)$ such that
$x_n=\E_n(x_\infty)$ for all $n \geq 1$. In this case, the sequence
$(x_n)_{n=1}^\infty$ converges to $x_\infty$ in $L^1(\M,\T)$. Similarly,
if $1<p<\infty$, every bounded $L^p$-martingale is of the form $(\E_n(x_\infty))_{n=1}^\infty$
for some $x_\infty \in L^p(\M,\T)$.

\medskip

The following decomposition of bounded $L^1$-martingale is the
non-commutative extension of the classical Krickeberg's
decomposition of martingales into linear combinations of positive
martingales. It will be used in the sequel. A proof for the finite
case can be found in \cite{Cuc} but the general case is readily
verified with the same techniques.
\begin{theorem}\label{decomposition} Let
$(x_{n})^{\infty}_{n=1}$ be a bounded $L^1$-martingale then $(x_{n})^{\infty}_{n=1}$
admits the following decomposition:
\begin{equation*}
x_{n} = \left(x^{(1)}_{n} - x^{(2)}_{n}\right) +
i\left(x^{(3)}_{n} - x^{(4)}_{n}\right)
\end{equation*}
for all $n\geq 1$  where for each $j\in \{1, 2, 3, 4\}$, the
sequence $(x^{(j)}_{n})^{\infty}_{n=1}$ is a  positive martingale.
Moreover, if $x_{n} = x^{*}_{n}$, for all $n\geq 1$,  then
$\sup_{n\geq 1}\|x_{n}\|_{1} = \T(x^{(1)}_{1}) + \T(x_{1}^{(2)})$.
\end{theorem}

We end this section  with a
maximal inequality type result.  Inspired by  Pisier's vector-valued non-commutative
$L^p$-spaces, Junge \cite{Ju} developed an abstract situation that can
efficiently describe a non-commutative analogue of the maximal function
theory
for bounded $L^p$-martingales when $p>1$. The proposition below can
be viewed as a substitute for the classical weak type $(1,1)$ boundedness
of maximal functions. Since it was not presented in the form below and
plays a crucial role in the proof of
 our main result,
  we will
reproduce the proof given in \cite{Cuc}.

\begin{proposition}\label{maximal}
If $(x_{n})^{\infty}_{n=1}$ is a positive bounded $L^1$-martingale
and $\lambda >0$ then there exists a sequence of decreasing projections
$(q_{n}^{(\lambda)})^{\infty}_{n=1}$ in $\M$ with:
\begin{itemize}
\item[(i)] for every $n \geq 1$, $q_{n}^{(\lambda)} \in \M_{n}$;
\item[(ii)] $q_{n}^{(\lambda)}$ commutes  with $q_{n-1}^{(\lambda)} x_{n}
q_{n-1}^{(\lambda)}$;
\item[(iii)] $q_{n}^{(\lambda)} x_{n} q_{n}^{(\lambda)} \leq \lambda
q_{n}^{(\lambda)}$;
\item[(iv)] $(q_{n}^{(\lambda)})_{n=1}^\infty$ is a decreasing sequence and
if we  set $q^{(\lambda)} = \bigwedge_{n=1}^{\infty}
q_{n}^{(\lambda)}$ then $\T ({\bf 1}-q^{(\lambda)}) \leq \T(x_{1})/\lambda$.
\end{itemize}
\end{proposition}

\begin{proof}
Let $q_{0}^{(\lambda)} = {\bf 1}$ and  inductively on $n \geq 1$,
define
\begin{equation*}
q_{n}^{(\lambda)} := \ch_{[0,\lambda]} (q_{n-1}^{(\lambda)} x_{n}
q_{n-1}^{(\lambda)}).
\end{equation*}
 The above definition makes
sense since $q_{n-1}^{(\lambda)} x_{n}q_{n-1}^{(\lambda)}$ is a
positive operator. It  is clear (by induction) that for every $n
\geq 1$, $q_{n}^{(\lambda)} \in \M_{n}$. Moreover, condition~(ii)
follows directly from the definition of $q_{n}^{(\lambda)}$ above.

For (iii), $q_{n}^{(\lambda)} x_{n} q_{n}^{(\lambda)} =
q_{n}^{(\lambda)} (q_{n-1}^{(\lambda)} x_{n} q_{n-1}^{(\lambda)})
q_{n}^{(\lambda)} = \ch_{[0,\lambda]} (q_{n-1}^{(\lambda)} x_{n}
q_{n-1}^{(\lambda)}) \cdot q_{n-1}^{(\lambda)} x_{n}
q_{n-1}^{(\lambda)} \leq \lambda q_{n}^{(\lambda)}$.  For (iv), it
is clear that $(q_{n}^{(\lambda)})_{n=1}^\infty$ is decreasing and
for every fixed  $n \geq 1$,
\begin{equation*}
\begin{split}
\T(x_{1}) &= \T(x_{n})\cr
&= \T(x_{n}q_{n}^{(\lambda)}) + \sum^{n}_{k=1} \T(x_{n}
(q_{k-1}^{(\lambda)}-q_{k}^{(\lambda)}))\cr
&= \T(q_{n}^{(\lambda)} x_{n} q_{n}^{(\lambda)}) + \sum^{n}_{k=1}
\T(\E_{k}
(x_{n})(q_{k-1}^{(\lambda)} - q_{k}^{(\lambda)})).
\end{split}
\end{equation*}
Since $\T(q_{n}^{(\lambda)} x_{n} q_{n}^{(\lambda)}) \geq 0$, we have
\begin{equation*}
\begin{split}
 \T(x_1) &\geq \sum^{n}_{k=1} \T ((q_{k-1}^{(\lambda)} - q_{k}^{(\lambda)}) x_{k}
(q_{k-1}^{(\lambda)} - q_{k}^{(\lambda)}))\cr &= \sum^{n}_{k=1}
\T((q_{k-1}^{(\lambda)} - q_{k}^{(\lambda)}) (q_{k-1}^{(\lambda)}
x_{k} q_{k-1}^{(\lambda)}) (q_{k-1}^{(\lambda)} -
q_{k}^{(\lambda)})).
\end{split}
\end{equation*}
From the definition of $q_{k}^{(\lambda)}$, it is clear that
$q_{k-1}^{(\lambda)} - q_{k}^{(\lambda)}= \ch_{(\lambda, \infty)}
(q_{k-1}^{(\lambda)} x_{k} q_{k-1}^{(\lambda)})$  and therefore
$(q_{k-1}^{(\lambda)} - q_{k}^{(\lambda)})q_{k-1}^{(\lambda)}
x_{k}q_{k-1}^{(\lambda)} (q_{k-1}^{(\lambda)} -
q_{k}^{(\lambda)})\geq \lambda (q_{k-1}^{(\lambda)} -
q_{k}^{(\lambda)})$ hence,
\begin{equation*}\T(x_1) \geq \lambda \sum^{n}_{k=1} \T
(q_{k-1}^{(\lambda)} - q_{k}^{(\lambda)})= \lambda \T({\bf
1}-q_{n}^{(\lambda)}).
\end{equation*}
Taking the limit as $n$ goes to $\infty$, (iv) follows. This
completes the proof.
\end{proof}

\section{Main Result}
In this section, we keep all notations introduced in the preliminaries.
In particular, all adapted sequences are understood to be with respect to a fixed
filtration of von Neumann subalgebras.
The following theorem answers positively a question raised by Pisier and Xu
\cite[Problem~7.5]{PX3} and is the main result of this paper.

\begin{theorem}\label{main}
There is an absolute constant $C$ such that if $x=(x_n)_{n=1}^\infty$ is a bounded
$L^1$-martingale and
$(\xi_{n})^{\infty}_{n=1}$
is an adapted sequence such that:
\begin{itemize}
\item[(i)] for every $n\geq 2$, $\xi_{n-1}$ commutes with  $\M_{n}$;
\item[(ii)] $\sup_{n\geq 1}\|\xi_{n}\|_\infty\leq 1$.
\end{itemize}
Then for every $N\geq 2$,
\begin{equation}\label{inequality1}
\left\Vert x_{1} + \sum^{N}_{k=2} \xi_{k-1} dx_k \right\Vert_{1,\infty}
\leq C
\left\Vert x_{N} \right\Vert_{1}.
\end{equation}
\end{theorem}

\medskip

\noindent
{\bf Proof of Theorem~\ref{main}:} By Theorem~\ref{decomposition}, it is enough to prove the case where
 $(x_{n})^{\infty}_{n=1}$ is a
positive martingale and $(\xi_{n})^{\infty}_{n=1}$ is an adapted
sequence of
self-adjoint operators satisfying conditions $(i)$ and $(ii)$. In the course of the proof, we will  frequently
use the tracial property of $\T$ and the $\T$-invariance property
of the expectations $\E_n$'s. For notational purpose, we set
$\xi_0={\bf 1}$.

Our goal is to show that there is a constant $C$, independent of
$(x_n)_{n=1}^\infty$ and $(\xi_n)_{n=1}^\infty$, such that for every
$\lambda>0$ and $N\geq 2$,
\begin{equation}\label{weaktype}
\T\left( \ch_{(\lambda, \infty)}\left(\left\vert \sum_{k=1}^N \xi_{k-1}dx_k
\right\vert\right)
\right) \leq \frac{C}{\lambda} \left\Vert x_N \right\Vert_1.
\end{equation}
The proof is divided into several steps:

\noindent {Step~1.} ({\it Reduction to bounded difference
sequences}). Fix $\lambda > 0$ and denote simply by
$(q_{n})^{\infty}_{n=1}$ (resp. $q$) the projections
$(q_{n}^{(\lambda)})^{\infty}_{n=1}$ (resp. $q^{(\lambda)})$ from
Proposition~\ref{maximal} and let $N \geq 2$ be fixed throughout
the proof.
\begin{lemma}\label{lemma1} For every $\alpha \in (0,1)$ and every
$\beta \in (0,1)$,
$$\T\left( \ch_{(\lambda, \infty)}\left(\left\vert
\sum_{k=1}^N \xi_{k-1}dx_k \right\vert\right)
\right) \leq
\alpha^{-1} \T\left( \ch_{(\beta\lambda, \infty)}\left(\left\vert \sum_{k=1}^N
q\xi_{k-1}dx_k q \right\vert\right)
\right) + \frac{2(1-\alpha)^{-1}}{\lambda}\T(x_1). $$
\end{lemma}
\begin{proof} We begin by splitting the operator
$S=\sum^{N}_{k=1} \xi_{k-1} dx_{k}$ into three  parts:
$$S=qSq + ({\bf 1}-q)Sq + S({\bf 1}-q).$$
 Fix $\alpha \in (0,1)$ and $\beta \in (0,1)$.
Using properties of the generalized singular value functions
$\mu(\cdot)$ from  \cite{FK},
\begin{equation*}
\begin{split}
\T \left(\ch_{(\lambda,\infty)} (|S|)\right) &= \int^{\infty}_{0}
\ch_{(\lambda,\infty)} \left\{\mu_{t}(S)\right\} \ dt \\
&\leq \int^{\infty}_{0} \ch_{(\lambda,\infty)} \left\{\mu_{\alpha
t}(qSq)
 + \mu_{(1-\alpha) t/2}(({\bf 1}-q)Sq)
+ \mu_{(1-\alpha) t/2} (S({\bf 1}-q)) \right\}\ dt \\
&= \int^{\infty}_{0} \ch_{(\lambda,\infty)} \left\{\mu_{\alpha
t}(qSq)
 + \mu_{(1-\alpha) t/2} (qS({\bf 1}-q))
+ \mu_{(1-\alpha) t/2}(S({\bf 1}-q)) \right\}  \ dt.
\end{split}
\end{equation*}
As $\mu_{(1-\alpha) t/2} (qS({\bf 1}-q)) \leq \mu_{(1-\alpha) t/2}
(|S({\bf 1}-q)|)$,
\begin{equation*}
\begin{split}
\T \left(\ch_{(\lambda,\infty)} (|S|)\right)
 &\leq \int^{\infty}_{0} \ch_{(\lambda,\infty)} \left\{\mu_{\alpha
t}(qSq)
 + 2\mu_{(1-\alpha) t/2} (|S({\bf 1}-q)|) \right\}  \ dt \\
&\leq \int^{\infty}_{0} \ch_{(\beta\lambda,\infty)}
\left\{\mu_{\alpha t} (qSq)\right\}\ dt
 + \int^{\infty}_{0} \ch_{((1-\beta)\lambda,\infty)}
\left\{\mu_{(1-\alpha) t/2}(2|S({\bf 1}-q)|)\right\} \ dt \\
&= \int^{\infty}_{0} \mu_{\alpha t}
\left\{\ch_{(\beta\lambda,\infty)} (|qSq|) \right\}\ dt +
 \int^{\infty}_{0} \mu_{(1-\alpha) t/2} \left\{\ch_{((1-\beta)\lambda,\infty)}
(2|S({\bf 1}-q)|)\right\}\ dt.
\end{split}
\end{equation*}
Remark that the projection $\ch_{((1-\beta)\lambda,\infty)} (2|S({\bf 1}-q)|)$
 is a  subprojection of $({\bf 1} -q)$ so
\begin{equation*}
\T \left(\ch_{(\lambda,\infty)} (|S|)\right) \leq
\int^{\infty}_{0} \mu_{\alpha t}
\left\{\ch_{(\beta\lambda,\infty)} (|qSq|) \right\}\ dt +
\int_{0}^{\infty} \mu_{(1-\alpha)t/2}({\bf 1}-q) \ dt
\end{equation*}
and by change of variables,
\begin{equation*}
\T \left(\ch_{(\lambda,\infty)} (|S|) \right)\leq
\alpha^{-1}\int^{\infty}_{0} \mu_{t}
\left\{\ch_{(\beta\lambda,\infty)} (|qSq|) \right\}\ dt +
2(1-\alpha)^{-1}\int_{0}^{\infty} \mu_{t}({\bf 1}-q) \ dt
\end{equation*}
which shows that
$\T \left(\ch_{(\lambda,\infty)} (|S|)\right) \leq
\alpha^{-1}\T \left(\ch_{(\beta\lambda,\infty)} (|qSq|)\right) +
2(1-\alpha)^{-1} \T(x_1)/\lambda.$
\end{proof}

\noindent Step~2. ({\it Reduction to difference sequence of a
supermartingale in $L^2(\M,\T)$}).
\begin{lemma}\label{lemma2}
The sequence $(q_kx_k q_k)_{k=1}^\infty$ is a supermartingale in
$L^2(\M,\T)$ and for every $\beta \in (0,1)$,
$$
\T\left(\ch_{(\beta\lambda,\infty)} \left(\left|\sum^{N}_{k=1} q \xi_{k-1}
dx_{k} q\right|\right)
\right) \leq
\frac{1}{\beta^2\lambda^{2}} \left\|q_{1} x_{1} q_{1} + \sum^{N}_{k=2} \xi_{k-1}
(q_{k} x_{k} q_{k} -
q_{k-1} x_{k-1} q_{k-1})  \right\|^{2}_{2}.
$$
\end{lemma}
\begin{proof} We remark first that since both sequences $(q_k)_{k=1}^\infty$
and $(x_k)_{k=1}^\infty$ are adapted, it is clear that
 $(q_{k} x_{k} q_{k})^{\infty}_{k=1}$ is adapted.
To prove  that it is
 a supermartingale, we need to verify that for every $k\geq 2$,
  $\E_{k-1} (q_{k}
x_{k} q_{k}) \leq q_{k-1} x_{k-1} q_{k-1}$. For this,  we remark
from Proposition~\ref{maximal} that
  since $q_{k}$ commutes with $q_{k-1} x_{k}q_{k-1}$
and $q_{k} \leq q_{k-1}$,
 $q_{k} x_{k} q_{k}
\leq q_{k-1} x_{k} q_{k-1}$. As $\E_{k-1}$ is a positive contraction,
\begin{equation*}
\begin{split}
\E_{k-1} (q_{k} x_{k} q_{k}) &\leq \E_{k-1} (q_{k-1} x_{k} q_{k-1})\cr
&= q_{k-1} \E_{k-1} (x_{k}) q_{k-1}\cr
&= q_{k-1} x_{k-1} q_{k-1}.
\end{split}\end{equation*}
For the second part of the lemma, it is clear that
 $$\T \left(\ch_{(\beta\lambda,\infty)} \left(\left|\sum^{N}_{k=1} q \xi_{k-1} dx_{k}
q\right|\right)\right) \leq \frac{1}{\beta^2\lambda^{2}} \T \left(\left|\sum^{N}_{k=1} q \xi_{k-1} dx_{k}
q\right|^{2}\right). $$  Moreover,
$
q\xi_{k-1} dx_{k} q = q \xi_{k-1} x_{k} q - q \xi_{k-1} x_{k-1} q
= q (q_{k} \xi_{k-1} x_{k} q_{k}) q - q (q_{k-1} \xi_{k-1} x_{k-1} q_{k-1})q
$.
Since $\xi_{k-1}$ commutes with $q_{k}$ and $q_{k-1}$, we conclude
that
$$ q\xi_{k-1} dx_{k} q=
 q \left(\xi_{k-1} (q_{k} x_{k} q_{k} - q_{k-1} x_{k-1} q_{k-1})\right)q.
$$
Similarly, $q dx_{1} q = q (q_{1} x_{1} q_{1}) q$ and therefore,
\begin{equation*}
\begin{split}
 \T \left(\left|\sum^{N}_{k=1}
q \xi_{k-1} dx_{k} q\right|^{2}\right) &=
   \left\|\sum^{N}_{k=1} q \xi_{k-1} dx_{k} q \right\|^{2}_{2} \\
&= \left\|q \left(q_{1} x_{1} q_{1} + \sum^{N}_{k=2} \xi_{k-1}
(q_{k} x_{k} q_{k} -
q_{k-1} x_{k-1} q_{k-1}) \right) q \right\|^{2}_{2} \\
&\leq\left\|q_{1} x_{1} q_{1} + \sum^{N}_{k=2} \xi_{k-1}
(q_{k} x_{k} q_{k} -
q_{k-1} x_{k-1} q_{k-1})  \right\|^{2}_{2}.
\end{split}\end{equation*}
This proves the lemma.
\end{proof}

\noindent
Step~3.({\it Change the supermartingale into sum of a martingale
and a decreasing sequence of operators}).
This is very standard: Define
\begin{equation}\label{Y}
 y_{k} :=\begin{cases} q_{1} x_{1} q_{1} \ &\text{for}\  k=1 \\
 q_{k} x_{k} q_{k} + \sum^{k-1}_{l=1} q_{l} x_{l} q_{l}
- \E_{l} (q_{l+1} x_{l+1} q_{l+1})\  &\text{for}\  k \geq 2.
\end{cases}
\end{equation}
Likewise,
\begin{equation}\label{Z}
z_{k} :=\begin{cases}  0 \ &\text{for}\ k=1 \\
\sum^{k-1}_{l=1} \E_{l} (q_{l+1} x_{l+1} q_{l+1}) -
q_{l} x_{l} q_{l} \ &\text{for} \ k \geq 2.
\end{cases}
\end{equation}
It is clear that $(y_k)_{k=1}^\infty$ is a positive martingale.
Moreover, for every $k\geq 1$,
\begin{equation}\label{YZ}
y_k + z_k = q_k x_k q_k
\end{equation}
and for every $k \geq 2$,
\begin{equation}\label{Z2}
z_{k} \leq z_{k-1} \leq \dots
 \leq z_{1} = 0.
\end{equation}
\begin{lemma}\label{lemma3}
The sequence $(y_{k})^{\infty}_{k=1}$ is a bounded
$L^{2}$-martingale with
$$\left\|y_{N}\right\|^{2}_{2} \leq
 6\lambda \T \left(q_{1} x_{N}\right) -
4\lambda \T \left(q_{N} x_{N}\right) - \left\|q_{1}
x_{1} q_{1}\right\|^{2}_{2} \leq 6 \lambda \T(x_1).$$
\end{lemma}
\begin{proof}
We will use the identity $ \left\|y_{N}\right\|^{2}_{2} =
\left\|(|y_{1}|^2 + \sum^{N}_{k=2}|y_{k} -
y_{k-1}|^{2})^{{1}/{2}}\right\|^{2}_{2} $. The main idea is to
estimate the sum $\sum_{k=2}^N\| y_k -y_{k-1}\|^{2}_2$ by a
telescopic sum. For $k\geq 2$, we notice  first from (\ref{Y})
that $y_{k} = y_{k-1} + q_{k} x_{k} q_{k} - \E_{k-1} (q_{k}
x_{k}q_{k})$ and therefore
\begin{equation*}
\begin{split}
y_{k} - y_{k-1} &= q_{k} x_{k} x_{k} - \E_{k-1} (q_{k} x_{k} q_k) \\
 &= \left(q_{k} x_{k} q_{k} - q_{k-1} x_{k-1} q_{k-1}\right) +
\left(q_{k-1} x_{k-1}q_{k-1} - \E_{k-1} (q_{k} x_{k} q_{k})\right).
\end{split}
\end{equation*}
Since $\|\cdot\|_{2}^2$ is convex,
\begin{equation*}
\begin{split}
\left\|y_{k} - y_{k-1}\right\|^{2}_{2} &\leq 2
(\left\|q_{k} x_{k} q_{k} - q_{k-1} x_{k-1} q_{k-1}\right\|^{2}_{2} +
 \left\|q_{k-1}
x_{k-1} q_{k-1} - \E_{k-1} (q_{k} x_{k} q_{k})\right\|^{2}_{2}) \\
 &= 2 \T \left((q_{k}
x_{k} q_{k} - q_{k-1} x_{k-1} q_{k-1})^{2}\right) +
2 \T \left((q_{k-1} x_{k-1}
q_{k-1} - \E_{k-1} (q_{k} x_{k} q_{k}))^{2}\right)\\
&= I + II.
\end{split}
\end{equation*}
We will estimate $I$ and $II$ separately. First for $I$, we use
the identity $(a-b)^{2} = a^{2} - b^{2} + b(b-a) + (b-a)b$ for
self-adjoint operators. With $a=q_kx_kq_k$ and
$b=q_{k-1}x_{k-1}q_{k-1}$, we have by taking the trace,
\begin{equation*}
\begin{split}
 I
&= 2\T\left((q_{k} x_{k} q_{k})^{2} -
(q_{k-1} x_{k-1} q_{k-1})^{2}\right) +
4 \T\left(q_{k-1} x_{k-1} q_{k-1}[q_{k-1} x_{k-1} q_{k-1} -
q_{k} x_{k} q_{k}]\right) \\
 &= 2 \T\left((q_{k} x_{k} q_{k})^{2} -
(q_{k-1} x_{k-1} q_{k-1})^{2}\right) +
4 \T\left(q_{k-1} x_{k-1} q_{k-1}[q_{k-1} x_{k-1} q_{k-1} -
\E_{k-1}(q_{k} x_{k} q_{k})]\right).
\end{split}
\end{equation*}
 By Proposition~\ref{maximal}~(iii),
$\left\|q_{k-1} x_{k-1} q_{k-1}\right\|_{\infty} \leq \lambda$.
Moreover, as $(q_kx_kq_k)_{k=1}^\infty$ is a supermartingale,
 $q_{k-1} x_{k-1} q_{k-1} - \E_{k-1} (q_{k} x_{k} q_{k}) \geq 0$.
 Therefore, we get
\begin{equation*}
\begin{split}
 I  &\leq
2\T\left((q_{k} x_{k} q_{k})^{2} -
(q_{k-1} x_{k-1} q_{k-1})^{2}\right) +
4 \lambda \T\left(q_{k-1} x_{k-1} q_{k-1} -
\E_{k-1}(q_{k} x_{k} q_{k})\right) \\
&=2\T\left((q_{k} x_{k} q_{k})^{2} -
(q_{k-1} x_{k-1} q_{k-1})^{2}\right) +
4 \lambda \T\left(q_{k-1} x_{k-1} q_{k-1} -
q_{k} x_{k} q_{k}\right).
\end{split}
\end{equation*}
For $II$, again since
 $q_{k-1} x_{k-1}
q_{k-1} \geq
q_{k-1} x_{k-1}
q_{k-1} - \E_{k-1} (q_{k} x_{k} q_{k}) \geq 0$, we have
 $$\left\|q_{k-1} x_{k-1}
q_{k-1} - \E_{k-1} (q_{k} x_{k} q_{k})\right\|_{\infty} \leq
\left\| q_{k-1} x_{k-1}
q_{k-1} \right\|_{\infty} \leq \lambda.$$
Hence, we get
\begin{equation*}
\begin{split}
  II &\leq
 2 \lambda \T \left(q_{k-1} x_{k-1} q_{k-1} -
\E_{k-1}(q_{k} x_{k} q_{k})\right) \\
&=  2 \lambda \T
\left(q_{k-1} x_{k-1} q_{k-1} - q_{k} x_{k} q_{k}\right).
\end{split}\end{equation*}

 Combining  the preceding estimates on $I$ and $II$, we conclude that
 for every  $k \geq 2$,
\begin{equation*}
\left\|y_{k} - y_{k-1}\right\|^{2}_{2}
\leq 2 \left(\left\|q_{k} x_{k} q_{k}\right\|^{2}_{2} -
\left\|q_{k-1} x_{k-1} q_{k-1}\right\|^{2}_{2}\right) +
6\lambda\T \left(q_{k-1} x_{k-1} q_{k-1} - q_{k} x_{k} q_{k}\right).
\end{equation*}
To conclude the proof of the lemma, we take the summation over $k$,
\begin{equation*}\begin{split}
\left\|y_{N}\right\|^{2}_{2}&=\left\|q_1 x_1 q_1\right\|^{2}_2 +
\sum_{k=2}^N \left\| y_k - y_{k-1}\right\|^{2}_2 \\
 &\leq \left\|q_{1} x_{1} q_{1}\right\|^{2}_{2} +
 2 \sum^{N}_{k=2}
(\left\|q_{k} x_{k} q_{k}\right\|^{2}_{2} -
\left\|q_{k-1} x_{k-1} q_{k-1}\right\|^{2}_{2}) \\
&\ +
6 \lambda \sum^{N}_{k=2} \T \left(q_{k-1} x_{k-1} q_{k-1} - q_{k} x_{k}
q_{k}\right) \\
&= \left\|q_{1} x_{1} q_{1}\right\|^{2}_{2} +
2 (\left\|q_{N} x_{N} q_{N}\right\|^{2}_{2} -
\left\|q_{1} x_{1} q_{1}\right\|^{2}_{2}) + 6\lambda
\T\left(q_{1} x_{1} q_{1} -
q_{N} x_{N} q_{N}\right) \\
&= 2\left\|q_{N} x_{N} q_{N}\right\|^{2}_{2} -
\left\|q_{1} x_{1} q_{1}\right\|^{2}_{2} +
6 \lambda \T \left((q_{1} - q_{N})x_{N}\right) \\
&\leq 2 \lambda \T\left(q_{N} x_{N}\right) - \left\|q_{1} x_{1}
q_{1}\right\|^{2}_{2} +
6\lambda \T \left((q_{1} - q_{N})x_{N} \right) \\
&= 6\lambda \T \left(q_{1} x_{N}\right) - 4 \lambda \T\left(q_{N}
x_{N}\right) - \left\|q_{1} x_{1} q_{1} \right\|^{2}_{2} \leq
6\lambda\T(x_1)
\end{split}\end{equation*}
which completes the proof.
\end{proof}

\noindent Step~4. ({\it Removal of the sequence
$(\xi_n)_{n=1}^\infty$ from the estimates}). This is done by
arguing separately on transforms of the difference sequences of
$(y_k)_{k=1}^\infty$ and $(z_k)_{k=1}^\infty$.
\begin{lemma}
$\left\|q_{1} x_{1} q_{1} + \sum^{N}_{k=2}
\xi_{k-1} (q_{k} x_{k} q_{k} - q_{k-1} x_{k-1} q_{k-1})\right\|^{2}_{2}
 \leq 4\left\|y_{N}\right\|_{2}^2$.
\end{lemma}
\begin{proof}
From the definitions of $(y_k)_{k=1}^\infty$ and $(z_k)_{k=1}^\infty$,
the convexity of $\Vert \cdot \Vert_{2}^2$ implies,
\begin{equation*}
\begin{split}
\left\|q_{1} x_{1} q_{1} + \sum^{N}_{k=2}
\xi_{k-1} (q_{k} x_{k} q_{k} - q_{k-1} x_{k-1} q_{k-1})\right\|^{2}_{2}
&\leq 2\left\| \sum^{N}_{k=1} \xi_{k-1} dy_k\right\|^{2}_{2} +
2\left\|\sum^{N}_{k=2} \xi_{k-1} (z_{k}-z_{k-1})\right\|^{2}_{2} \\
&= III +  IV.
\end{split}
\end{equation*}
As in Step~3, we will estimate $III$ and $IV$ separately. First,
since martingale transforms are clearly bounded (with constant=1) in
$L^{2}(\M,\T)$, it follows that
\begin{equation*}
III \leq 2\left\| \sum^{N}_{k=1} dy_k\right\|^{2}_{2} =2\left\|y_N
\right\|^{2}_2
\end{equation*}
which gives an upper bound of $III$ that is independent of the
sequence $(\xi_k)_{k=1}^\infty$.

 On the other hand, it is clear that
\begin{equation*}
\begin{split}
IV &=
2\left\|\sum^{N}_{k=2} \xi_{k-1} (z_{k}-z_{k-1})\right\|^{2}_2 \\
&=2 \T\left(\left|\sum^{N}_{k=2} \xi_{k-1}
(z_{k-1}-z_{k})\right|^{2} \right) \\
&= 2\T\left(\sum^{N}_{k=2}\sum_{l=2}^N
(z_{k-1}-z_{k}) \xi_{k-1} \xi_{l-1} (z_{l-1} -z_{l})\right) \\
&=2  \sum^{N}_{k=2}\sum_{l=2}^N \T\left(
(z_{k-1}-z_{k}) \xi_{k-1} \xi_{l-1} (z_{l-1} -z_{l})\right).
\end{split}
\end{equation*}
To estimate $IV$, recall from (\ref{Z2})
 that $z_{k-1} - z_{k}$ and $z_{l-1} - z_{l}$ are positive operators.
Assume for instance that $k \leq l$ (the case $l \leq k$ is
handled equally) then by assumption, $\xi_{k-1} \xi_{l-1}$
commutes with $\M_{k}$ so we have,
$$ (z_{k-1} - z_{k}) \xi_{k-1} \xi_{l-1}= (z_{k-1} -
z_{k})^{1/2} \xi_{k-1} \xi_{l-1} (z_{k-1} -z_{k})^{1/2}
\leq z_{k-1} - z_{k}.$$
Therefore by taking the trace,
\begin{equation*}
\begin{split}
\T \left((z_{k-1} - z_{k}) \xi_{k-1} \xi_{l-1} (z_{l-1} -
z_{l})\right)  &= \T\left((z_{l-1} - z_{l})^{1/2}
[(z_{k-1}-z_k)\xi_{k-1} \xi_{l-1}] (z_{l-1} -z_{l})^{1/2})
\right) \\
&\leq \T\left((z_{l-1} -
z_{l})^{1/2} (z_{k-1}-z_k) (z_{l-1} -z_{l})^{1/2})\right) \\
&=
 \T\left((z_{k-1}-z_{k})(z_{l-1}-z_{l})\right).
\end{split}
\end{equation*}
Hence, we get
\begin{equation*}
\begin{split}
IV &= 2  \sum^{N}_{k=2}\sum_{l=2}^N \T\left(
(z_{k-1}-z_{k}) \xi_{k-1} \xi_{l-1} (z_{l-1} -z_{l})\right) \\
&\leq 2  \sum^{N}_{k=2}\sum_{l=2}^N
\T\left((z_{k-1}-z_{k})(z_{l-1} -z_{l})\right) \\
&= 2\T\left((\sum^{N}_{k=2} z_{k}-z_{k-1})^{2}\right) \\
&= 2\Vert z_N \Vert_{2}^2.
\end{split}
\end{equation*}

By combining the preceding estimates on $III$ and $IV$, we obtain
\begin{equation*}
\left\|q_{1} x_{1} q_{1} + \sum^{N}_{k=2}
\xi_{k-1} (q_{k} x_{k} q_{k} - q_{k-1} q_{k-1} x_{k-1})\right\|^{2}_{2}
 \leq 2\left\|y_{N}\right\|_{2}^2 + 2\left\|z_{N}\right\|^{2}_{2}.
\end{equation*}
To conclude the proof of the lemma, note from (\ref{YZ}) that
$y_{N} - q_{N} x_{N} q_{N} = -z_{N} \geq 0$ so $y_{N} \geq -z_{N}
\geq 0$ which implies $\left\|y_{N}\right\|_{2} \geq
\left\|z_{N}\right\|_{2}$.
\end{proof}

To complete the proof of Theorem~\ref{main}, it is  enough, as
mentioned above, to verify (\ref{weaktype}). This is obtained by
putting together the four lemmas above. Indeed,
\begin{equation*}
\begin{split}
\T\left( \ch_{(\lambda, \infty)}\left(\left\vert
\sum_{k=1}^N \xi_{k-1}dx_k \right\vert\right)
\right) &\leq
\alpha^{-1} \T\left( \ch_{(\beta\lambda, \infty)}\left(\left\vert \sum_{k=1}^N
q\xi_{k-1}dx_k q \right\vert\right)
\right) + \frac{2(1-\alpha)}{\lambda}^{-1}\T(x_1) \\
&\leq \frac{\alpha^{-1}}{\beta^2\lambda^{2}} \left\|q_{1} x_{1} q_{1} +
\sum^{N}_{k=2} \xi_{k-1}
(q_{k} x_{k} q_{k} -
q_{k-1} x_{k-1} q_{k-1})  \right\|^{2}_{2} \\
&\ +
\frac{2(1-\alpha)^{-1}}{\lambda}\T(x_1) \\
&\leq  \frac{4\alpha^{-1}\beta^{-2}}{\lambda^{2}}\left\|y_{N}\right\|_{2}^2 +
\frac{2(1-\alpha)^{-1}}{\lambda}\T(x_1) \\
&\leq \frac{24\alpha^{-1}\beta^{-2} + 2(1-\alpha)^{-1}}{\lambda}\T(x_1).
\end{split}
\end{equation*}
This shows that (\ref{weaktype}) is satisfied with
$$C=\inf\left\{24\alpha^{-1}\beta^{-2} + 2(1-\alpha)^{-1};\
\alpha \in (0,1),\ \beta \in (0,1)\right\}=
 \frac{14\sqrt{3}}{3} +28.$$
\qed

\begin{remark}
In the proof above, no significant  effort was made to minimize
the constant $C$ involved in Theorem~\ref{main}.  Recall that in
the classical case, the sharp constant $C=2$ is known and was
obtained by Burkholder in \cite{Bu2}.  His approach, as expected,
is based on a stopping time argument which (at least at the time
of this writing) does not seem to have an efficient
non-commutative analogue.
\end{remark}

\begin{problem}
 Find the {\lq \lq{sharp}"} constant $C$ for which
 Inequality~(\ref{inequality1}) holds?
\end{problem}

Theorem~\ref{main} can be extended to transforms of submartingales
and supermartingales .
\begin{corollary}
There exists a constant $K$ such that if  $(x_n)_{n=1}^\infty$ is
either a submartingale or a supermartingale and is  bounded in
$L^1(\M,\T)$ then for any sequence of signs
$(\epsilon_n)_{n=1}^\infty$,
\begin{equation*}
\sup_{N\geq 2}\left\|\epsilon_1x_{1} + \sum^{N}_{k=2}
\epsilon_{k}(x_k- x_{k-1})\right\|_{1,\infty} \leq K \sup_{n\geq
1} \left\|x_{n}\right\|_{1}.
\end{equation*}
\end{corollary}
\begin{proof}
We will present the proof for submartingale. As in the proof of
Theorem~\ref{main}, we split $(x_n)_{n=1}^\infty$ into sum of a
martingale and an increasing sequence of positive operators. Let
\begin{equation*}
 y_{k} :=\begin{cases}  x_{1}  \ &\text{for}\  k=1 \\
  x_{k}  + \sum^{k-1}_{l=1}  x_{l}
- \E_{l} ( x_{l+1} )\  &\text{for}\  k \geq 2.
\end{cases}
\end{equation*}
and
\begin{equation*}
z_{k} :=\begin{cases}  0 \ &\text{for}\ k=1 \\
\sum^{k-1}_{l=1} \E_{l} ( x_{l+1} ) -  x_{l}  \ &\text{for} \ k
\geq 2.
\end{cases}
\end{equation*}
The following properties are immediate:
\begin{itemize}
\item[(a)] $(y_k)_{k=1}^\infty$ is a  martingale;
\item[(b)] for every $k \geq 1$,
$y_k + z_k =  x_k $;
\item[(c)]
 for every $k \geq 2$,
$z_{k} \geq z_{k-1} \geq \dots
 \geq z_{1} = 0$.
\end{itemize}
Moreover, for every $k\geq 1$,
\begin{equation*}
\begin{split}
\Vert z_k \Vert_1 &= \T(z_k) \\
&=\sum_{l=1}^{k-1} \T( \E_l(x_{l+1}) -x_l)\\
&=\sum_{l=1}^{k-1} \T( x_{l+1} -x_l)\\
&=\T(x_{k-1}-x_1) \leq 2\Vert x_k \Vert_1.
\end{split}
\end{equation*}
As above,
\begin{equation*}
\begin{split}
\left\|\epsilon_1x_{1} + \sum^{N}_{k=2} \epsilon_{k}(x_k-
x_{k-1})\right\|_{1,\infty} &\leq 2\left\| \sum^{N}_{k=1}
\epsilon_{k}dy_k\right\|_{1,\infty} +2 \left\| \sum^{N}_{k=2}
\epsilon_{k}(z_k- z_{k-1})\right\|_{1,\infty} \\
&\leq 2C\left\| \sum^{N}_{k=1} dy_k\right\|_{1} + 2 \left\|
\sum^{N}_{k=2} \epsilon_{k}(z_k- z_{k-1})\right\|_{1}.
\end{split}
\end{equation*}
It is easy to see that $-z_N \leq \sum^{N}_{k=2} \epsilon_{k}(z_k-
z_{k-1}) \leq z_N$. Therefore $\Vert \sum^{N}_{k=2}
\epsilon_{k}(z_k- z_{k-1})\Vert_1 \leq \Vert z_N \Vert_1$ and
hence
\begin{equation*}
\begin{split}
\left\|\epsilon_1x_{1} + \sum^{N}_{k=2} \epsilon_{k}(x_k-
x_{k-1})\right\|_{1,\infty} &\leq 2C\Vert y_N \Vert_1 + 2\Vert
z_N\Vert_1 \\
&\leq 2C\Vert x_N\Vert_1 + (2C +2)\Vert z_N \Vert_1 \leq K\Vert
x_N \Vert_1.
\end{split}
\end{equation*}
The proof is complete.
\end{proof}

 As in the
commutative case, Theorem~\ref{main} implies that if $\T({\bf
1})<\infty$,  martingale transforms are bounded from $L^1(\M,\T)$
into $L^p(\M,\T)$ for $0<p<1$.

\begin{corollary} Assume that $\T({\bf 1})<\infty$. Under the assumption of Theorem~\ref{main},
for every $0<p<1$, there exists a constant $K_{p}$ (depending only on $p$)
such that:
\begin{equation*}
\left\|x_{1} + \sum^{N}_{k=2} \xi_{k-1}dx_k\right\|_{p} \leq K_{p}
\left\|x_{N}\right\|_{1}.
\end{equation*}
\end{corollary}

In \cite{PX}, Pisier and Xu proved, as a consequence of the
non-commutative Burkholder-Gundy inequalities, a non-commutative
analogue of Stein's inequality (\cite[Theorem~2.3]{PX},
\cite[Theorem~8 p.~103]{St}) for $1<p<\infty$.  Their proof
reveals that what is needed is the unconditionality of martingale
differences in $L^{p}(\M,\T)$.  A slightly different proof was
given by Junge and Xu (\cite{JX}) which yields a better constant.
Below, we will adopt their proof together with Theorem~\ref{main}
to get the corresponding result for $p=1$.

\begin{theorem}\label{stein}
There is a constant $\gamma >0$ such that for any
finite sequence $(a_{k})^{n}_{k=1}$ in
$L^1(\M,\T)$,
\begin{equation*}
\left\|\left(\sum^{n}_{k=1} \E_{k} (a_{k})^{*} \E_{k}
(a_{k})\right)^{{1}/{2}}\right\|_{1,\infty} \leq \gamma
\left\|\left(\sum^{n}_{k=1} a_{k}^{*}
a_{k}\right)^{{1}/{2}}\right\|_{1}.
\end{equation*}
\end{theorem}

\begin{proof}
Consider the tensor product $(\M,\T) \otimes (B(\ell^{2}_{n}), \sigma)$
where $\sigma = n^{-1}tr$ is the usual normalized trace on $B(\ell^{2}_{n})$.
For $k\geq 1$, let $\widetilde{\E_{k}} = \E_{k} \otimes Id_{B(\ell^{2}_{n})}$ be the
conditional expectation from  $\M \otimes B(\ell^{2}_{n})$ onto the
subalgebra $\M_{k} \otimes B(\ell^{2}_{n})$.

Let $A_{k} = n a_{k} \otimes e_{k,1}$ for $1\leq k \leq n$ and
$(r_{j})_{j \geq
1}$ be the sequence of the Rademacher functions on $[0,1]$.  Then for any
$t \in [0,1]$,
\begin{equation*}
\sum^{n}_{k=1} \E_{k} (a_{k})^{*} \E_{k} (a_{k}) \otimes
n^2e_{1,1} = \left|\sum^{n}_{k=1}
\widetilde{\E}_{k}(r_{k}(t)A_{k})\right|^{2}
\end{equation*}
and therefore
\begin{equation*}
\begin{split}
\left\|\left(\sum^{n}_{k=1} \E_{k}(a_{k})^{*} \E_{k}
(a_{k})\right)^{{1}/{2}}\right\|_{1,\infty} &= \left\|\sum^{n}_{k=1}
\widetilde{\E}_{k}(r_{k}(t)A_{k})\right\|_{1,\infty}\cr
&= \left\|\sum^{n}_{k=1} \widetilde{\E}_{n}(r_{k}(t)A_{k}) - \sum^{n-1}_{k=1}
\sum^{n-1}_{j=1}(\widetilde{\E}_{j-1} -
\widetilde{\E}_{j})(r_{k}(t)A_{k})\right\|_{1,\infty}.
\end{split}
\end{equation*}
Since $\|a +b \|_{1,\infty} \leq 2\|a\|_{1,\infty} +
2\|b\|_{1,\infty}$ for every $a$ and $b$ in $L^{1,\infty}(\M,\T)$,
\begin{equation*}
\begin{split}
\left\|\left(\sum^{n}_{k=1} \E_{k}(a_{k})^{*} \E_{k}
(a_{k})\right)^{{1}/{2}}\right\|_{1,\infty}
&\leq 2\left\|\sum^{n}_{k=1} \widetilde{\E}_{n}(r_{k}(t)A_{k})\right\|_{1,\infty} +
2\left\|\sum^{n-1}_{k=1} \sum^{n-1}_{j=1} (\widetilde{\E}_{j-1} -
\widetilde{\E}_{j})(r_{k}(t)A_{k})\right\|_{1,\infty}\cr
&\leq 2\left\|\sum^{n}_{k=1} \widetilde{\E}_{n}(r_{k}(t)A_{k})\right\|_{1} +
2\left\|\sum^{n-1}_{k=1} \sum^{n-1}_{j=1}(\widetilde{\E}_{j-1} -
\widetilde{\E}_{j})(r_{k}(t)A_{k})\right\|_{1,\infty}\cr
&\leq 2\left\|\sum^{n}_{k=1} r_{k}(t)A_{k} \right\|_{1} +2\left\|\sum^{n-1}_{k=1}
\sum^{n-1}_{j=1}(\widetilde{\E}_{j-1} - \widetilde{\E}_{j})(r_{k}(t)A_{k})
\right\|_{1,\infty}.
\end{split}
\end{equation*}
Let $f = \sum^{n-1}_{k=1} r_{k} A_{k}$ and consider the filtration $(\M_{k}
\otimes B(\ell^{2}_{n}) \otimes L^\infty(\cal{F}_{k}))_{k \geq 1}$ where
$\cal{F}_{k}$ is the
$\sigma$-field generated by $\{r_{1},r_{2},\dots,r_{k}\}$.  Denoting by
$(df_{j})_{j \geq 1}$ the difference sequence of $f$ with respect to this
filtration, we have:
$$
\sum^{n-1}_{j=1} (\widetilde{\E}_{j-1} - \widetilde{\E}_{j})
(\sum^{j}_{k=1} r_{k} (t) A_{k}) = \sum^{n-1}_{j=1} df_{2j+1}.
$$
By Theorem~\ref{main}, we conclude that
\begin{equation*}
\begin{split}
\left\|\left(\sum^{n}_{k=1} \E_{k} (a_{k})^{*} \E_{k}
(a_{k})\right)^{{1}/{2}}\right\|_{1,\infty} &\leq 2 \left\|\left(\sum^{n}_{k=1} a_{k}^{*}
a_{k}\right)^{{1}/{2}}\right\|_{1} + 2C\left\|f \right\|_{1}\cr
&\leq (2 + 2C)\left\|\left(\sum^{n}_{k=1} a^{*}_{k} a_{k}\right)^{{1}/{2}}\right\|_{1}.
\end{split}
\end{equation*}
This shows  the theorem with $\gamma \leq 2+2C$.
\end{proof}

\section{Estimating UMD-constants for non-commutative  spaces}

In this section, we are primarily interested in UMD-constants of
non-commutative spaces. Our main motivation comes mainly from a
question of Pisier \cite{PS6} on the order  of the UMD-constants
of the Schatten class $S^p$. To this end, we began by reviewing
the relevant background on UMD-spaces (UMD stands for
unconditional martingale differences).
\begin{definition}
A Banach space $X$ is said to have the UMD-property if for some $p
\in (1, \infty)$, there exists a constant $C$, which depends only
on $p$ and $X$ such that for all $n\geq 1$,
\begin{equation}\label{UMD}  \left\Vert \sum_{j=1}^n \epsilon_j d_j
\right\Vert_{L^p(X)} \leq C \left\Vert \sum_{j=1}^n  d_j
\right\Vert_{L^p(X)}
\end{equation}
for every $X$-valued martingale
difference sequence $(d_j)_{j=1}^\infty$ and
$(\epsilon_j)_{j=1}^\infty \in \{-1,1\}^{\N}$.
\end{definition}
Here $L^p(X)= L^p(\Omega, \Sigma, \mu,  X)$ denotes the Bochner
space of all strongly measurable functions $f$ on a probability
space $(\Omega, \Sigma,\mu)$ with values in $X$ such that:
$$\Vert f
\Vert_{L^p(X)} :=(\int_\Omega \Vert f(\omega)\Vert_{X}^p\
d\mu(\omega))^{1/p} <\infty.$$

Since we are interested in estimating the constants involved, we
will make distinctions between the indices. We will denote the
best constant in (\ref{UMD}) by $C_p(X)$. By duality, it is clear
that $X$ is a UMD-space if and only if $X^*$ is a UMD-space. In
this case, $C_p(X)= C_q(X^*)$ with $1/p + 1/q =1$. For more
information on UMD-spaces, we refer to \cite{BO10} and \cite{Bu3}.
\begin{theorem}{\rm (\cite{PS6})}
Let $X$ be a UMD-space then for any $1<p,q<\infty$, there exist
positive constants $\alpha(p,q)$ and $\beta(p,q)$ depending only
on $p$ and $q$ such that:
\begin{equation*}
\alpha(p,q)C_p(X) \leq C_q(X) \leq  \beta(p,q) C_p(X).
\end{equation*}
In particular,  for any $p\geq 3$, we have $(2\sqrt{3})^{-1}
C_2(X) \leq C_p(X) \leq 7pC_2(X)$.
\end{theorem}

Our main tool in this section is  the unconditionality of
martingale transforms on $L^{p}(\M,\T)$ (for $ 1 < p < \infty$)
which follows from our main result.  More precisely,
\begin{theorem}\label{transformp}
Let $1 < p < \infty$.  For any finite non-commutative
$L^{p}$-martingale $x$ and any sequence of signs
$(\epsilon_n)_{n=1}^\infty$,
$$
\left\|\sum_{n \geq 1} \epsilon_{n} dx_{n}\right\|_{p} \leq
c_{p}\|x\|_{p},$$ where $c_{p} \leq C{p^{2}}/(p-1)$ with $C$ being
a universal constant.
\end{theorem}

The case $1 < p < 2$ follows by interpolation from
Theorem~\ref{main} and the $L^{2}$-boundedness of martingale
transforms.  The case $2 < p < \infty$ can be deduced by duality.
\begin{remark}
Except for the constants, Theorem~\ref{transformp} was obtained in
\cite{PX}.  As $c_{p} \leq Cp^{2}/(p-1)$, it is clear that
$c_{p}=O(p)$ when $p \to \infty$ and $O((p-1)^{-1})$ when $p \to
1$. These are the optimal order of growths for $c_{p}$.
\end{remark}
 We will
apply Theorem~\ref{transformp} to estimate the  UMD-constants of $L^p(\M,\T)$.
It is well known that for $1<p<\infty$, $L^p(\M,\T)$ (and in
particular the Schatten class $S^p$) is a UMD-space. This was
established as a consequence of the characterization of UMD-spaces
due to Burkholder \cite{Bu3} and Bourgain \cite{BO10} in
terms of vector-valued Hilbert transforms (\cite{BGM},
\cite{BO10}). Such approach gives constants that are $O(p^2)$ when
$p \to \infty$. We remark that the UMD property of $L^p(\M,\T)$ also follows
from the generalized Riesz projections associated with group representations
which was extensively studied by Zsid\'o \cite{ZS}.
Our next result follows immediately from
Theorem~\ref{transformp} and the definition of UMD-spaces. It
answers positively a question from \cite{PS6}.

\begin{corollary}\label{UMDconstant}
There exists a constant $C$  such that for every $1<p<\infty$,
\begin{equation*}
C_p(L^p(\M,\T)) \leq Cp^2/(p-1).
\end{equation*}
In particular, there exists a constant $C'$ such that for $p\geq
2$, $C_p(L^p(\M,\T)) \leq C'p$.
\end{corollary}
\begin{proof}
Let $(\Omega, \Sigma,\mu)$ be a probability space and
$(d_n)_{n=1}^\infty$  be a $p$-integrable  $L^p(\M,\T)$-valued
martingale difference sequence defined on $(\Omega, \Sigma, \mu)$
relative to an increasing sequence of $\sigma$-subalgebras
$(\Sigma_n)_{n=1}^\infty$ of $\Sigma$ with conditional
expectations $(\mathbb{E}_n)_{n=1}^\infty$. Set
$\cal{N}=L^\infty(\Omega,\Sigma,\mu) \overline{\otimes} \M$ and
let $\cal{N}_n = L^\infty(\Omega,\Sigma_n,\mu) \overline{\otimes}
\M$. Then the conditional expectation $\E_n$ from $\cal{N}$ onto
$\cal{N}_n$ is given by $\mathbb{E}_n \otimes Id$. It is clear
that $(d_n)_{n=1}^\infty$ is a non-commutative martingale
difference sequence in $L^p(\cal{N}, \mu \otimes \T)$ associated
to the filtration $(\cal{N}_n)_{n=1}^\infty$. It is well known that
$L^p(\cal{N}, \mu \otimes \T)$ is isometrically isomorphic to
the Bochner space $L^p(\mu, L^p(\M,\T))$. By
Theorem~\ref{transformp}, for every $k\geq 1$ and $\epsilon_n=\pm
1$,
\begin{equation*}
\begin{split}
\left\Vert \sum_{n=1}^k \epsilon_n d_n \right\Vert_{L^p(\mu,
L^p(\M,\T))}&=\left\Vert \sum_{n=1}^k
\epsilon_n d_n \right\Vert_{L^p(\cal{N}, \mu \otimes \T)}\\
 &\leq c_p \left\Vert \sum_{n=1}^k
d_n \right\Vert_{L^p(\cal{N}, \mu \otimes \T)}
\end{split}
\end{equation*}
 which shows that $C_p(L^p(\M,\T)) \leq c_p$.
\end{proof}

\begin{remarks}
(1) The preceding corollary shows in particular that $C_p(S^p)$ and
$C_2(S^p)$ are $O(p)$ when $p \to \infty$.

(2) Replacing $L^\infty(\Omega, \Sigma, \mu)$ by a general non-commutative probability
space (in the sense of \cite[p.~48]{PIS7}), the proof of
Corollary~\ref{UMDconstant}
shows that the constant for the operator space version of UMD ($UMD_p$ property,
\cite[Definition~4.8]{PIS7}) of $L^p(\M,\T)$ is also $O(p)$ when $p \to \infty$.

(3) The constants relative to the boundedness of the $L^p(\M,\T)$-valued Hilbert transforms
are also $O(p)$ when $p \to \infty$ but this fact seems to provide only  weaker estimates
that  $C_2(L^p(\M,\T))$ is $O(p^2)$ when $p \to \infty$.
\end{remarks}

The above result can be extended to the Haagerup $L^p$-spaces
associated to general von Neumann algebras (we refer to \cite{HAA},
\cite{TE} for in depth description of such spaces) modulo the
following approximation of the Haagerup $L^p$-spaces.
\begin{theorem}{\rm (\cite{HAA2})}
Let $\M$ be an arbitrary von Neumann algebra and $L^p(\M)$ be the
Haagerup $L^p$-space associated with $\M$ ($0<p<\infty$). There
exist a Banach space $X$ (a $p$-Banach space if $0<p<1$), a
directed family $\left\{(\M_i, \T_i)\right\}_{i \in I}$ of finite
von Neumann algebras $\M_i$ (with normal faithful finite traces
$\T_i$), and a family $\left\{j_i\right\}_{i \in I} $ of isometric
embeddings $j_i: L^p(\M_i,\T_i) \to X$ such that:
\begin{itemize}
\item[(i)] $j_i\left( L^p(\M_i, \T_i)\right) \subset
j_k\left( L^p(\M_k, \T_k)\right)$ for all $i, k \in I$ with $i\leq
k$;
\item[(ii)] $\bigcup_{i\in I}j_i\left( L^p(\M_i, \T_i)\right)$ is
dense in $X$;
\item[(iii)] $L^p(\M)$ is isometric to a (complemented for $1\leq
p <\infty$) subspace of $X$.
\end{itemize}
\end{theorem}
Let $\M$ be an arbitrary von Neumann algebra (not necessarily
semi-finite) and $p>1$. If $X$ is the Banach  space obtained from
the above theorem then $X$ is a UMD-space with $C_p(X)=\sup_{i \in
I} C_p(L^p(\M_i, \T_i))$. In particular, the Haagerup $L^p$-space
$L^p(\M)$ is a UMD-space with constants equal to those of the
finite case.

Let us now consider the case $p=1$.  If $p=1$ or
$p=\infty$ then $S^p$ fails the UMD-property. Let us denote by
$S^1(n \times \infty)$ (resp. $S^1(\infty \times n)$) the space of
trace class operators for $n \times \infty$ matrices (resp.
$\infty \times n$ matrices). The next result gives an estimate of
the UMD-constant of $S^1(n \times \infty)$ when $ n \to \infty$.
It should be compared with \cite[Theorem~6.1]{PS6}.
\begin{theorem}
There exists a constant $K$ such that for any $n\geq 1$, we have
\begin{equation*}
C_2(S^1(\infty \times n)) \leq K\log(n+1)
\end{equation*}
and similarly for $S^1(n \times \infty)$.
\end{theorem}
\begin{proof}
For every $x \in S^1(\infty \times n)$ and $q\leq 2\leq p$ with
$1/q +1/p =1$,
\begin{equation*}
\left\Vert x \right\Vert_q \leq \left\Vert x \right\Vert_1 \leq
n^{1/p}\left\Vert x \right\Vert_q.
\end{equation*}
Hence
\begin{equation*}
C_2(S^1(\infty \times n)) \leq n^{1/p}C_2(S^q)
\end{equation*}
but since $C_2(S^q)=C_2(S^p) \leq  2\sqrt{3} C_p(S^p) \leq
2\sqrt{3}C'p$,
\begin{equation*}
C_2(S^1(\infty \times n)) \leq 2\sqrt{3}C'pn^{1/p}.
\end{equation*}
Choosing $p=\max\{2,\log(n)\}$, the theorem follows.
\end{proof}
We remark that since $S^1(n \times \infty)$ is the dual of the
space of operators $B(\ell^2_n,\ell^2)$ then
$C_2(B(\ell^2_n,\ell^2))$ is of order $\log(n)$. In particular, if
$M_{n,m}$ is the space of $n \times m$ matrices with the usual
norm then $C_2(M_{n,m})$ is of order $\min\{\log(n),\log(m)\}$.
The preceding argument also shows that for $N\geq 1$, there exist
a constant $K>0$ such that if $(x_n)_n$ is a finite martingale in
$S^1_N$ (as predual of $M_N$) then
$$\left\Vert\sum_n \epsilon_n dx_n \right\Vert_1 \leq K\log(N+1)\sup_n\left\Vert x_n \right\Vert_1$$
 for all $\epsilon_n=\pm 1$.

\medskip

We end this section by considering the general case of rearrangement invariant Banach function spaces.
Before proceeding, we need to recall the notion of Boyd indices.
Let $E$ be a rearrangement invariant Banach space on $(0,\infty)$.
For $s>0$, the dilation operator $D_s: E \to E$ is defined by setting
\begin{equation*}
D_sf(t)=f(t/s), \qquad t>0, \qquad f \in E.
\end{equation*}
The {\it lower and upper Boyd indices } of $E$ are defined by
\begin{equation*}
\underline{\alpha}_E :=\lim_{s\to 0^{+}}\frac{\log\Vert D_s\Vert}{\log s}, \qquad
\overline{\alpha}_E :=\lim_{s\to \infty}\frac{\log\Vert D_s\Vert}{\log s}.
\end{equation*}
It is well known that $0\leq \underline{\alpha}_E \leq \overline{\alpha}_E \leq 1$ and if
$E=L^p$  for $1\leq p \leq \infty$ then $\underline{\alpha}_E =\overline{\alpha}_E=1/p$.
If $0< \underline{\alpha}_E \leq \overline{\alpha}_E < 1$, we shall say that $E$ has
non-trivial Boyd indices. The next result can be viewed as a martingale analogue of \cite[Theorem~4.1]{DDPS}
where the existence of generalized Riesz projections where considered.
\begin{theorem}\label{transformE}
Let $E$ be a rearrangement invariant Banach function space on $(0,\infty)$ with Fatou norm.
The following statements are equivalent.
\begin{itemize}
\item[(i)] $E$ has non-trivial Boyd indices;
\item[(ii)] There exists a constant $c(E)$ depending only on $E$ such that for any  semi-finite von Neumann
algebra $(\M,\T)$ and any martingale $(x_n)_{n=1}^\infty$ in $E(\M,\T)$,
\begin{equation*}
\left\Vert \sum_{n=1}^N \epsilon_n dx_n \right\Vert_{E(\M,\T)} \leq c(E)
\left\Vert \sum_{n=1}^N dx_n \right\Vert_{E(\M,\T)}
\end{equation*}
for every $N\geq 2$ and $\epsilon_n=\pm 1$.
\end{itemize}
\end{theorem}
\begin{proof}
$(i) \implies (ii)$. Choose $1<p<q <\infty$ such that $1/q <\underline{\alpha}_E \leq \overline{\alpha}_E<1/p$
then $E$ is an interpolation space of the pair $(L^p, L^q)$ and therefore $E(\M,\T)$ is an interpolation space for the pair
$(L^p(\M,\T), L^q(\M,\T))$. The implication $(i) \implies (ii)$ follows by interpolation from Theorem~\ref{transformp}.

$(ii) \implies (i)$. Assume first that $\overline{\alpha}_E=1$. Choose a filtration of $B(\ell^2)$,  a finite martingale
$(x_n)_{n=1}^J$ in $S^1$
and a sequence $\epsilon_n=\pm 1$, $1\leq n \leq J$ such that
\begin{equation*}
\Vert x_J \Vert_1 =1, \qquad  \left\Vert \sum_{n=1}^J \epsilon_n dx_n \right\Vert_1 \geq 2c(E).
\end{equation*}
We can assume that $(x_n)_{n=1}^J$ are $N\times N$-matrices. Since $\overline{\alpha}_E=1$, it follows from
\cite[Proposition~2.b.6]{LT}, that there exist non-negative, disjointly supported, equidistributed functions
$(f_i)_{i=1}^N$ with $\Vert f_i \Vert_E=1 $ for $1\leq i \leq N$, and
\begin{equation*}
\frac{2}{3}\sum_{i=1}^N |a_i| \leq \left\Vert \sum_{i=1}^N a_if_i \right\Vert_E
\end{equation*}
for every choice of scalars $(a_i)_{i=1}^N$ in $\mathbb{C}$. Let $\M$ be $L^\infty(0,\infty) \otimes M_N(\mathbb{C})$ with
the trace $\T$ given by $\lambda \otimes tr$ where $\lambda$ denotes the trace on $L^\infty(0,\infty)$ induced by the Lebesgue
measure and $tr$ is the canonical trace on $M_N(\mathbb{C})$. We observe that
\begin{equation*}
\left\Vert f_1 \otimes A \right\Vert_{E(\M,\T)} = \left\Vert \sum_{i=1}^N s_i(A) f_i \right\Vert_E
\end{equation*}
for any matrix $A$ in $M_N(\mathbb{C})$, where $(s_i(A))_{i=1}^N$ denotes the singular values of $A$ arranged in decreasing order.
In fact, let $A$ be any $N\times N$ matrix and consider $D$ the diagonal matrix with entries $s_1(A), s_2(A), \dots, s_N(A)$.
If $U$ and $V$ are unitary matrices for which $A=UDV$ then for every $t>0$,
\begin{equation*}
\mu_t(f_1\otimes A)=\mu_t(1 \otimes U. f_1 \otimes D. 1 \otimes V)=\mu_t(f_1\otimes D)=
\mu_t\left(\sum_{i=1}^N s_i(A)f_j\right),
\end{equation*}
where the last equality follows from the fact that $(f_i)_{i=1}^N$ are equidistributed and the definition of the trace on $\M$, and this
proves the assertion.
It now follows that
\begin{equation*}
\Vert f_1 \otimes x_J \Vert_{E(\M,\T)} = \left\Vert \sum_{i=1}^N s_i(x_J)f_i\right\Vert_E \leq
\left(\sum_{i=1}^N s_i(x_J)\right)\Vert f_1\Vert_E=1.
\end{equation*}
On the other hand,
\begin{equation*}
\begin{split}
\left \Vert f_1 \otimes \sum_{n=1}^J \epsilon_n dx_n \right\Vert_{E(\M,\T)} &=
\left\Vert \sum_{i=1}^N s_i\left(\sum_{n=1}^J \epsilon_n dx_n \right)f_i \right\Vert_E \\
&\geq \frac{2}{3}\sum_{i=1}^N s_i\left(\sum_{n=1}^J \epsilon_n dx_n \right) \\
&=\frac{2}{3}\left\Vert \sum_{n=1}^J \epsilon_n dx_n \right\Vert_1 \geq \frac{4}{3} c(E).
\end{split}
\end{equation*}
Observe that $(f_1 \otimes x_n)_{n=1}^J$ is a finite martingale in $E(\M,\T)$. Assertion $(iii)$ implies that
\begin{equation*}
\left \Vert f_1 \otimes \sum_{n=1}^J \epsilon_n dx_n \right\Vert_{E(\M,\T)} \leq c(E) \Vert f_1 \otimes x_J \Vert_{E(\M,\T)},
\end{equation*}
and this yields a contradiction.
The same argument can be applied to prove that assertion $(iii)$ implies that $\underline{\alpha}_E >0$.
\end{proof}

Unlike the case of $L^p(\M,\T)$, Theorem~\ref{transformE} does not lead to UMD-property for $E(\M,\T)$. Special characterizations
that provide ready recognition of  UMD-property for rearrangement invariant Banach function spaces on $(0,\infty)$ seem to be unavailable.
On the other hand, there are examples of separable rearrangement invariant spaces on $(0,\infty)$ with non-trivial Boyd indices
which are not reflexive (see for instance \cite[p.~132]{LT}), and therefore fail the UMD-property. It is still an open  question
if $E$ being a UMD-space is sufficient for  $E(\M,\T)$ to be  a UMD-space.

\section{Non-commutative Burkholder-Gundy inequalities revisited}

In this section, we will point out that the weak-type inequality in
our main result
 implies the non-commutative Burkholder-Gundy inequalities proved in
\cite{PX}.  We first recall the two square functions introduced in
\cite{PX}.

Fix $1 \leq p < \infty$ and let $x$ be a bounded $ L^{p}$-martingale.
  Recall,
$$
S_{C,n} (x) = \left(\sum^{n}_{k=1}|dx_{k}|^{2}\right)^{{1}/{2}}
\ \ \text{and}\ \ \
S_{R,n}(x) = \left(\sum^{n}_{k=1}|dx_{k}^{*}|^{2}\right)^{{1}/{2}}.
$$
For any finite sequence  $a=(a_n)_{n \geq 1}$  in $L^{p}(\M,\T)$, set
$$
\|a\|_{L^{p}(\M;l^{2}_{C})} = \left\|\left(\sum_{n \geq
1}|a_{n}|^{2}\right)^{{1}/{2}}\right\|_{p}, \ \
\|a\|_{L^{p}(\M; l^{2}_{R})} = \left\| \left(\sum_{n \geq 1}
|a_{n}^{*}|^{2}\right)^{{1}/{2}}\right\|_{p}.
$$
The difference sequence $dx$ belongs to $L^{p}(\M;l^{2}_{C})$
(resp. $L^{p}(\M; l^{2}_{R})$)  if and only if  the sequence
$(S_{C,n}(x))_{n=1}^\infty$ (resp. $(S_{R,n}(x))_{n=1}^\infty$) is
a bounded  in $L^{p}(\M,\T)$.  In this case,  the limits
$S_{C}(x) = (\sum^{\infty}_{k=1}|dx_{k}|^{2})^{{1}/{2}}$ and
$S_{R}(x) = (\sum^{\infty}_{k=1}|dx_{k}^{*}|^{2})^{{1}/{2}}$ are
elements of $L^{p}(\M,\T)$.

For $1 \leq p < \infty$,  $\H^{p}_{C}(\M)$
 (resp. $\H^{p}_{R}(\M)$) is defined as
the set of all $L^{p}$-martingales $x$ with respect to $(\M_{n})_{n \geq
1}$ such that $dx \in L^{p}(\M; l^{2}_{C})$
(resp. $L^{p}(\M;l^{2}_{R})$),
and set
$$
\|x\|_{\H^{p}_{C}(\M)} = \|dx\|_{L^{p}(\M;l^{2}_{C})} \ \ \text{and}
\ \ \|x\|_{\H^{p}_{R}(\M)} = \|dx\|_{L^{p}(\M;l^{2}_{R})}.
$$
Equipped with the previous norms, $\H^{p}_{C}(\M)$ and $\H^{p}_{R}(\M)$
are Banach spaces.
The Hardy space of non-commutative martingale is defined as follows:  if $1
\leq p < 2$,
$$\H^{p}(\M) = \H^{p}_{C}(\M) + \H^{p}_{R}(\M)$$
 equipped with
the norm
$$
\|x\|_{\H^{p}(\M)} =
\inf\left\{\|y\|_{\H^{p}_{C}(\M)} + \|z\|_{\H^{p}_{R}(\M)}
: \ x = y + z, \ y \in \H^{p}_{C}(\M),\  z \in \H^{p}_{R}(\M) \right\};
$$
and if $2 \leq p < \infty $,
$$
\H^{p}(\M) = \H^{p}_{C}(\M) \cap \H^{p}_{R}(\M)$$
 equipped with the norm
$$\|x\|_{\H^{p}(\M)} = \max\left\{
\|x\|_{\H^{p}_{C}(\M)} , \ \|x\|_{\H^{p}_{R}(\M)}\right\}.$$
The main result
of \cite{PX} states that:
\begin{theorem}\label{BG}
Let $1 < p < \infty$.  Let $x = (x_{n})^{\infty}_{n=1}$ be an
$L^{p}$-martingale.  Then $x$ is
bounded in $L^{p}(\M,\T)$ if and only if $x$ belongs to $\H^{p}(\M)$.  If
this is the case then
\begin{equation*}
 \alpha^{-1}_{p}\|x\|_{\H^{p}(\M)} \leq \|x\|_{p} \leq
\beta_{p}\|x\|_{\H^{p}(\M)}. \tag {$BG_p$}
\end{equation*}
\end{theorem}

The strategy of \cite{PX2} and \cite{PX} for the particular cases
of tensor products,
Clifford algebras and the Free group von Neumann algebras was to show the
unconditionality of martingale differences in $L^{p}(\M,\T)$ (for $1<p<\infty$)
 using transference
argument to change  non-commutative martingales into  commutative
vector-valued ones, and then apply non-commutative Khintchine
inequalities (which we will recall below) together with a
non-commutative analogue of Stein's inequality. Such approach highlights the
fact that non-commutative $L^p$-spaces are UMD-spaces. Their proof for
the general case was completely different as they argued
inductively on $p=2^{n}$ for $n\geq 1$, then used interpolations
and duality.

\medskip

Let us recall the non-commutative Khintchine inequalities for the convenience
of the reader. Let $\epsilon=(\epsilon_n)_{n\geq 1}$ be a sequence of
independent random variables on some probability space $(\Omega, \cal{F},  P)$
such that $P(\epsilon_n=1)=P(\epsilon_n =-1)=1/2$ for all $n\geq 1$.
\begin{theorem}\label{Kintchine}
{\rm (Non-commutative Khintchine inequalities, \cite{LP4, LPI})}
Let $1 \leq p < \infty$.  Let $a = (a_{n})_{n \geq 1}$ be a finite
sequence in $L^{p}(\M,\T)$.
\begin{itemize}
\item[(i)] If $2 \leq p < \infty$,
\begin{equation*}
\begin{split}
\|a\|_{L^{p}(\M;l^{2}_{C}) \cap L^{p} (\M;l^{2}_{R})} &\leq
\left(\int_{\Omega}\|\sum_{n \geq 1} \epsilon_{n} a_{n}\|^{2}_{p}
\ dP(\epsilon)\right)^{{1}/{2}}\cr
&\leq \beta \sqrt{p} \|a\|_{L^{p}(\M;l^{2}_{C}) \cap
L^{p}(\M;l^{2}_{R}).}
\end{split}
\end{equation*}
\item[(ii)] If $1 \leq p < 2$,
\begin{equation*}
\begin{split}
\alpha \|a\|_{L^{p}(\M;l^{2}_{C}) + L^{p} (\M;l^{2}_{R})} &\leq
\left(\int_{\Omega}\|\sum_{n \geq 1} \epsilon_{n} a_{n}\|^{2}_{p}
\ dP(\epsilon)\right)^{{1}/{2}}\cr
&\leq \|a\|_{L^{p}(\M;l^{2}_{C}) + L^{p}(\M;l^{2}_{R})},
\end{split}
\end{equation*}
\end{itemize}
where $\alpha >0$ and $\beta >0$ are absolute constants.
\end{theorem}

As in the case of unconditionality of martingale difference
sequences, the non-commutative Stein's inequality can also be
deduced from Theorem~\ref{stein} above and interpolation. This
approach produces better estimate of the constant involved.

\begin{theorem}\label{steinp}
Let $1 < p < \infty$.  Define the map $Q$ on all finite sequences
$a=(a_{n})_{n \geq 1}$ in $L^{p}(\M,\T)$ by $Q(a)=(\E_{n}(a_{n}))_{n \geq
1}$.  Then
\begin{equation*}
\|Q(a)\|_{L^{p}(\M;l^{2}_{C})} \leq
\gamma_{p}\|a\|_{L^{p}(\M;l^{2}_{C})},\ \quad
\|Q(a)\|_{L^{p}(\M;l^{2}_{R})}
\leq \gamma_{p}\|a\|_{L^{p}(\M;l^{2}_{R})}, \tag {$S_p$}
\end{equation*}
where $\gamma_{p} \leq K{p^{2}/(p-1)}$ for some absolute constant $K$.
\end{theorem}
As noted in \cite{PX}, Theorem~\ref{steinp} shows that
$Q$ extends to a bounded linear projection on $L^{p}(\M;l^{2}_{C})$ and
$L^{p}(\M;l^{2}_{R})$. Consequently,
$\H^p(\M)$ is complemented in $L^{p}(\M;l^{2}_{C}) + L^{p}(\M;l^{2}_{R})$ or
$L^{p}(\M;l^{2}_{C}) \cap L^{p}(\M;l^{2}_{R})$
 according to $1 < p \leq 2$ or $2 \leq p <\infty$.

\medskip

We are now ready to present the proof.

\noindent
{\bf Proof of Theorem~\ref{BG}:}
  Let $1 < p < 2$. By Theorem~\ref{transformp} and Theorem~\ref{Kintchine},
\begin{equation*}
\alpha\|dx\|_{L^{p}(\M;l^{2}_{C}) + L^{p}(\M;l^{2}_{R})} \leq
c_{p}\|x\|_{p}.
\end{equation*}
Applying Theorem~\ref{transformp} to the martingale difference
$(\epsilon_{n} dx_{n})_{n=1}^\infty$
instead of $(dx_{n})_{n=1}^\infty$, we also have  the converse inequality:
\begin{equation*}
\|x\|_{p} \leq c_{p}\left\|\sum_{n \geq 1} \epsilon_{n} dx_{n}\right\|_{p}
\end{equation*}
for all
$\epsilon_{n} = \pm1$.
By Theorem~\ref{Kintchine},
\begin{equation*}
\|x\|_{p} \leq c_{p}\|dx\|_{L^{p}(\M;l^{2}_{C}) +
L^{p}(\M;l^{2}_{R})}
\end{equation*}
and therefore
\begin{equation*}
\alpha c^{-1}_{p}\|dx\|_{L^{p}(\M;l^{2}_{C}) +
L^{p}(\M;l^{2}_{R})} \leq \|x\|_{p} \leq
c_{p}\|dx\|_{L^{p}(\M;l^{2}_{C}) + L^{p}(\M;l^{2}_{R})}.
\end{equation*}
By duality, if  $2 < p$, then
\begin{equation*}
 c^{-1}_{p}\|dx\|_{L^{p}(\M;l^{2}_{C}) \cap L^{p}(\M;l^{2}_{R})}
\leq \|x\|_{p} \leq \alpha^{-1}c_{p}\|dx\|_{L^{p}(\M;l^{2}_{C})
\cap L^{p}(\M;l^{2}_{R})}.
\end{equation*}
This shows $(BG_{p})$ for $2 < p < \infty$ with $\alpha_{p} \leq
c_{p}$ and $\beta_{p} \leq \alpha^{-1}c_{p}$.

For the case $1 < p < 2$, remark  that $\|x\|_{\H^{p}(\M)} \geq
\|dx\|_{L^{p}(\M;l^{2}_{C}) + L^{p}(\M;l^{2}_{R})}$.  From
$(S_{p})$, we conclude that
\begin{equation*}
(\gamma_p)^{-1} c_{p}^{-1}\|x\|_{\H^{p}(\M)} \leq \|x\|_{p} \leq
c_{p}\|x\|_{\H^{p}(\M)}.
\end{equation*}
This proves $(BG_{p})$ for $1 < p < 2$ with $\alpha_{p} \leq
\gamma_pc_p$ and $\beta_{p} \leq c_{p}$. \qed

\begin{remarks}\label{growth}
As $\gamma_p \leq Kp^2/(p-1)$, $\gamma_p=O(p)$ when $p\to \infty$
and $O((p-1)^{-1})$ when $p\to 1$. These are the optimal orders
for $\gamma_p$.
 Recall that in the commutative case, the optimal order of
growths for the constants $\alpha_{p}$ and $\beta_{p}$ are (see
for instance \cite{Bu1}): $\beta_{p}$ is bounded when $p \to 1$
and $O(p)$ when $p \to \infty$; $\alpha_{p}$ is $O((p-1)^{-1})$
when $p \to 1$ and $O(\sqrt{p})$ when $p \to \infty$. The fact
that $\beta_{p}$ is bounded when $p \to 1$ for the non-commutative
case was recovered by Junge and Xu \cite[Corollary~4.3]{JX}.
Pisier showed in \cite{PS5} that $\beta_{p}$ is $O(p)$ for $p$
even integers.  The proof above also gives $\beta_{p}$ is $O(p)$
when $p \to \infty$. As for $\alpha_{p}$, the preceding proof
gives $\alpha_{p}$ is $O((p-1)^{-2})$ when $p \to 1$ and $O(p)$
when $p \to \infty$. For more in depth discussion about the orders
of growth of these constants, we refer to a recent paper of Junge
and Xu \cite{JX2}.
\end{remarks}

\section{Remarks on the class $L\log{L}$}

Recall first the class $L\log{L}$. If $L^0(\Omega,\cal{F}, P)$ is the space
of all (classes) of measurable functions on a given probability space
$(\Omega, \cal{F}, P)$, the class $L\log{L}$ is defined by setting
$$L\log{L}=\left\{ f \in L^0(\Omega, \cal{F}, P);
\int |f|\log^+|f|\ dP <\infty\right\}.$$ Set $\|f\|_{L\log{L}}
=\int |f|\log^+|f|\ dP$. Equipped with the equivalent norm $\| f
\| =\int_{0}^1 f^*(t) \log(1/t)\ dt$,
 the space $L\log{L}$ is a rearrangement invariant Banach function space
(see for instance \cite[Theorem~6.4, pp.~246-247]{BENSHA}) so a
non-commutative analogue $L\log{L}(\M,\T)$ is well defined as
described in Sect.~2. We remark that if a martingale $x$ is
bounded in $L\log{L}(\M,\T)$ then it is uniformly integrable in
$L^1(\M,\T)$ and therefore is of the form
$x=(\E_n(x_\infty))_{n=1}^\infty$ with $x_\infty \in
L\log{L}(\M,\T)$.

The starting point of this section is the following well known inequality
from the classical  theory.
\begin{theorem}
There is a constant $K$ such that  if
 $(f_{k})^{\infty}_{k=1}$ is a (commutative) martingale  then
for every $n \geq 1$,
\begin{equation}\label{log}
\mathbb{E}\left(\sup_{1 \leq k \leq n}|f_{k}|\right) \leq K +
K\mathbb{E}\left(|f_{n}|\log^{+}|f_{n}|\right).
\end{equation}
\end{theorem}
By the equivalence of maximal functions and square functions for
(commutative) martingales \cite{Da2}, the left hand side of
(\ref{log}) can be replaced by $\mathbb{E}(S_n(f))$ where
$S_n(f)=(\sum^{n}_{k=1}|df_{k}|^{2})^{{1}/{2}}$.  The standard
procedure for establishing inequality~(\ref{log}) above is to
derive first the weak-type inequality for maximal functions by a
stopping time argument then integrating from $1$ to $\infty$ (see
\cite[pp.~317-318]{Doob};  consult also \cite[pp.~81-85]{GA}
 for another approach). In a more operator  theoretical point of view, inequality~(\ref{log})
follows from general theory of interpolation of operators of weak
types (see for instance \cite[Theorem~6.6, pp.~248-249]{BENSHA}).
With this observation, the following result follows immediately
from Theorem~\ref{main}:
\begin{theorem}\label{LlogL1}
There exists a constant $K$ such that if $x=(x_{n})^{\infty}_{n=1}$ is
a martingale which is bounded in $L\log{L}(\M,\T)$ then for any
sequence of signs $(\epsilon_{n})^{\infty}_{n=1}$,
\begin{equation*}
\sup_{N\geq 1}\left\|\sum^{N}_{n=1} \epsilon_{n} dx_{n}\right\|_{1}
\leq K  + K\left\| x_\infty  \right\|_{L\log{L}(\M,\T)}.
\end{equation*}
\end{theorem}
Using the non-commutative Khintchine inequality, one can deduce
\begin{corollary}\label{LlogL2}
There is a constant $K$ such that if $x=(x_{n})^{\infty}_{n=1}$ is a
martingale that is bounded in $L\log{L}(\M,\T)$, then
\begin{equation*}
\left\|dx\right\|_{L^{1}(\M;l^{2}_{C}) + L^{1}(\M;l^{2}_{R})} \leq K +
K\left\| x_\infty  \right\|_{L\log{L}(\M,\T)}.
\end{equation*}
\end{corollary}

Corollary~\ref{LlogL2} can be viewed as a non-commutative
extension of (\ref{log}) above. However, inequality (\ref{log}) is
equivalent to:
 if $f$ is bounded in $L\log{L}$ then $f \in \H^{1}$.
Since $\left\| x\right\|_{\cal{H}^1(\M)} \geq
\left\|dx\right\|_{L^{1}(\M;l^{2}_{C}) + L^{1}(\M;l^{2}_{R})}$,
the following question arises naturally:

\begin{problem}
  Does there exist a constant
$K$ such that for every martingale $x$:
$$\left\|x\right\|_{\H^{1}(\M)} \leq K
+K\left\|x_\infty\right\|_{L\log{L}(\M,\T)}?$$
\end{problem}

An old argument from conjugate function theory together with the
fact noted in Remark~\ref{growth} above that $\alpha_p$ is
$O((p-1)^{-2})$ when $p\to 1$ can be used to prove a related
inequality. The proof given below is modelled after a presentation
in Zygmund's book (\cite[p.119]{ZYG}).
\begin{proposition} There is an absolute constant $K$ such that if
$x=(x_n)_{n=1}^\infty$ is a martingale that is bounded in
$L\log{L}$ and $\T(|x_\infty| (\log^+|x_\infty|)^2)<\infty$, then
\begin{equation*}
\Vert x \Vert_{\cal{H}^1(\M)} \leq K + K\T\left(|x_\infty|
(\log^+|x_\infty|)^2\right).
\end{equation*}
\end{proposition}
\begin{proof}
Let $x=(\E_n(x_\infty))_{n=1}^\infty$ be a martingale  with
 $\T(|x_\infty| (\log^+|x_\infty|)^2)<\infty$.
 Let $a = |x_\infty|$
 and  set $(e_t)_t $ to be the spectral
decomposition of $a$. For each $k \in \N$, let $P_k =\chi_{[2^{k-1},2^k)}(a)$
be the  spectral projection relative to $[2^{k-1}, 2^k)$. Define $a_k =aP_k$
 for
$k \geq 1$ and $a_0=a\chi_{[0, 1)}(a)$.
 Clearly
$a= \sum_{k=0}^{\infty} a_k $ in $L^1(\M,\T)$.

  For every $k\in\N$, consider the martingale
$x^{(k)}=(\E_n(x_\infty P_k))_{n=1}^\infty$ then
 $\|x^{(k)}\|_{\cal{H}^1(\M)} \leq
\|x^{(k)}\|_{\cal{H}^p(\M)} \leq \alpha_p\|x^{(k)}\|_p$. So for
every $1<p<2$, there is a constant $C$ such that,
$\|x^{(k)}\|_{\cal{H}^1(\M)} \leq
C^2p^4(p-1)^{-2}\|x^{(k)}\|_{p}$. Since $\|x^{(k)}\|_{p} =\|
a_k\|_p$ and
  $a_k \leq 2^k P_k$, we get for $1<p<2$,
$$\|x^{(k)}\|_{\cal{H}^1(\M)} \leq
  16C^2 (p-1)^{-2} 2^k \T(P_k)^{\frac{1}{p}}.$$
If we set $p =1 + 1/(k+1)$ and $\eta_k =\T(P_k)$, we have
 $$ \|x^{(k)}\|_{\cal{H}^1(\M)}  \leq 16C^2 (k+1)^2 2^k \eta_{k}^{\frac{k+1}{k+2}}.$$
Taking the summation over $k$,
 $$\|x\|_{\cal{H}^1(\M)}  \leq \sum_{k=0}^\infty 16C^2(k+1)^2 2^k \eta_{k}^{\frac{k+1}{k+2}}.$$
We note as in \cite{ZYG} that if $J=\{ k \in \N;\ \eta_k \leq 3^{-k}\}$
then
 $$\sum_{k \in J} 16C^2(k+1)^2 2^k \eta_{k}^{\frac{k+1}{k+2}} \leq
 \sum_{k=0}^\infty 16C^2(k+1)^2 2^k (3^{-k})^{\frac{k+1}{k+2}} = \alpha <\infty.$$
On  the other hand, for $k \in \N \setminus J$,
$\eta_{k}^{\frac{k+1}{k+2}} \leq \eta_k 3^{\frac{k}{k+2}} \leq
\beta \eta_k$ where $\beta=\sup_k 3^{\frac{k}{k+2}}$. So we get
\begin{align*}
\|x\|_{\cal{H}^1(\M)}  &\leq
 \alpha + 16C^2\beta \sum_{k=0}^\infty (k+1)^2 2^k \eta_k \\
 &\leq \alpha +16C^2\beta(\eta_0 + 8\eta_1) +16C^2\beta \sum_{k\geq 2}(k+1)^2
2^k \eta_k.
\end{align*}
Since for $k\geq 2$,\   $k+1 \leq 3(k-1)$, we get
$$\|x\|_{\cal{H}^1(\M)} \leq
\alpha +128C^2\beta +288C^2\beta\sum_{k\geq2}(k-1)^2 2^{k-1}\eta_k.$$
To complete the proof, notice that  for $k\geq 2$,
\begin{align*}
 (k-1)^2 2^{k-1} \eta_k &= \int_{2^{k-1}}^{2^k}  (k-1)^2 2^{k-1}\ d\T(e_t)\\
  &\leq \int_{2^{k-1}}^{2^k} \frac{t(\log t)^2}{(\log 2)^2} \ d\T(e_t),
\end{align*}
as  $2^{k-1} \leq t $ and therefore $(k-1)\log 2 \leq \log t$.
Hence if we set
$$K=\max\{\alpha +128C^2\beta, 288C^2\beta(\log 2)^{-2}\},$$
then  we get:
 $$\|x\|_{\cal{H}^1(\M)} \leq  K + K \T\left(a (\log^+(a))^2\right).$$
The proof is complete.
\end{proof}

We remark that combining Corollary~\ref{LlogL2} and Theorem~\ref{stein},
one can deduce the following:
There exists a constant $K'$ such that:
\begin{equation*}
\begin{split}
&\inf\left\{\|dy\|_{L^{1,\infty}(\M;l^{2}_C)} +
\|dz\|_{L^{1,\infty}(\M;l^{2}_R)}
: \ x = y + z, \ y \in \H^{1}_{C}(\M),\  z \in \H^{1}_{R}(\M) \right\} \\
&{}\quad \leq K' + K'\left\|x_\infty\right\|_{L\log{L}(\M,\T)}.
\end{split}
\end{equation*}
The next question corresponds to the weak type boundedness of
square functions:

\begin{problem}
  Does there exist a constant $K$ such that for every
bounded $L^1$-martingale $x$,
$$
\inf\left\{\|dy\|_{L^{1,\infty}(\M;l^{2}_C)} +
\|dz\|_{L^{1,\infty}(\M;l^{2}_R)}
: \ x = y + z, \ y \in \H^{1}_{C}(\M),\  z \in \H^{1}_{R}(\M) \right\}
\leq  K\left\|x \right\|_1 ?
$$
\end{problem}

We remark that a simple adjustment of  the proof of Theorem~\ref{main}
gives: There exists a constant $K$ such that
for every $\lambda < 0$,
\begin{equation*}
\inf\{\lambda \T(\ch_{(\lambda,\infty)}(S_{C}(y))) +
\lambda \T(\ch_{(\lambda,\infty)}
(S_{R}(z))):\  x=y+z \} \leq K\|x\|_{1}.
\end{equation*}
We conclude by noticing that the proof of Theorem~\ref{stein} combined
with Theorem~\ref{LlogL1} yields the following:
\begin{theorem} There exists a constant $K$ such that
for any  finite sequence $a=(a_k)^{n}_{k=1}$ in $L\log{L}(\M,\T)$, if
$Q(a)=(\E_{k}(a_{k}))^{n}_{k=1}$ then
\begin{equation*}
\|Q(a)\|_{L^{1}(\M;l^{2}_{C})} \leq \|a\|_{L^{1}(\M; l^{2}_{C})} +
K +
K \T\left((\sum^{n}_{k=1}|a_{k}|^{2})^{{1}/{2}}
\log^{+}(\sum^{n}_{k=1}|a_{k}|^{2})^{{1}/{2}}\right).
\end{equation*}
\end{theorem}

\medskip

\noindent {\bf Acknowledgements.} I am very grateful to Q.~Xu for
several comments on an earlier version of this paper.

\providecommand{\bysame}{\leavevmode\hbox
to3em{\hrulefill}\thinspace}
\providecommand{\MR}{\relax\ifhmode\unskip\space\fi MR }
\providecommand{\MRhref}[2]{%
  \href{http://www.ams.org/mathscinet-getitem?mr=#1}{#2}
} \providecommand{\href}[2]{#2}

\end{document}